\documentclass{jctart08}
\usepackage{graphics, graphicx}     %%%   Пакеты для работы с рисунками.

\usepackage{cite}
\usepackage{amsmath}
\usepackage{amssymb}
\usepackage{amsthm}
\theoremstyle{plain}
\newtheorem{Theorem}{Theorem}[section] %
\newtheorem{Lemma}{Lemma}[section]

\theoremstyle{definition}
\renewcommand{\thetable}{\thesection.\arabic{table}}

\addto\captionsrussian{\def\refname{References}}
\newenvironment{Proof} % имя окружения
{\par\noindent{\it Proof of}} % команды для \begin
{\hfill$\vspace{5mm}\scriptstyle\blacksquare$} % команды для \end

\numberwithin{equation}{section} %Чтобы нумерация в каждой секции была независимой
\numberwithin{figure}{section} %Чтобы нумерация в каждой секции была независимой
\numberwithin{table}{section} %Чтобы нумерация в каждой секции была независимой

\begin{document}

\setcounter{page}{1}

\markboth{M.I. Isaev}{Asymptotic behaviour of the number of Eulerian circuits}

\title{Asymptotic behaviour of the number of Eulerian circuits}
\date{}
\author{ {\bf M.I. Isaev}\\
Moscow Institute of Physics and Technology,
141700 Dolgoprudny, Russia\\
Centre de Math\'ematiques Appliqu\'ees, Ecole Polytechnique,
91128 Palaiseau, France\\
e-mail: \tt{isaev.m.i@gmail.com}}

\maketitle
{\bf Abstract}
\begin{abstract}
	We determine the asymptotic behaviour of the number of Eulerian circuits in undirected simple 
	graphs with large algebraic connectivity (the second-smallest eigenvalue of the Laplacian matrix). 
	We also prove some new properties of the Laplacian matrix.

\end{abstract}

\section{Introduction}
$\ \ \ $

Let $G$ be a simple connected graph all of whose vertices have even degree. A {\it Eulerian
circuit } in $G$ is a closed walk (see, for example, \cite{Biggs1976}) which uses every edge of $G$ exactly once. We let
$Eul(G)$ denote the number of these up to cyclic equivalence. Our purpose in this paper is
to estimate $Eul(G)$ for those $G$ having large algebraic connectivity.

Our method is to adopt the proof given in \cite{Brendan1995} for the case $G = K_n$. We refer to that paper
for the interesting history of this problem, and suggest that readers who want to understand
our proofs carefully may find it helpful to have a copy at hand. Since the publication of
\cite{Brendan1995}, the work \cite{Brightwell2004} has appeared showing that counting the number of Eulerian circuits in an
undirected graph is complete for the class $\#P$. Thus this problem is difficult in terms of
complexity theory.

	Here is an outline of the paper. The asymptotic
formula for $Eul(K_n)$ and our main result are presented and discussed in Section 2. 
In Section 3 we prove some basic properties of the Laplacian matrix, which may be of independent interest.  
In Section 4 we express $Eul(G)$ in terms of an
$n$-dimensional integral using Cauchy’s formula. 
The value of the integral is estimated in Sections 5 and 6, using some Lemmas proved in Section 8. 
We prove the main result in Section 7.

%%%%%%%%%%%%%%%%%%%%%%%%%%%%%%%%%%%%%%%%%%%%%%%%%%%%%%%%%%%%%%%%%%%%%%%%%%%%%%%%%%%%%%%%%%%%%%%%%%%%%%%%%%%%%%%%%%%%%%%%%
%%%%%%%%%%%%%%%%%%%%%%%%%%%%%%%%%%%%%%%%%%%%%%%%%%%%%%%%%%%%%%%%%%%%%%%%%%%%%%%%%%%%%%%%%%%%%%%%%%%%%%%%%%%%%%%%%%%%%%%%%
%%%%%%%%%%%%%%%%%%%%%%%%%%%%%%%%%%%%%%%%%%%%%%%%%%%%%%%%%%%%%%%%%%%%%%%%%%%%%%%%%%%%%%%%%%%%%%%%%%%%%%%%%%%%%%%%%%%%%%%%%
%%%%%%%%%%%%%%%%%%%%%%%%%%%%%%%%%%%%%%%%%%%%%%%%%%%%%%%%%%%%%%%%%%%%%%%%%%%%%%%%%%%%%%%%%%%%%%%%%%%%%%%%%%%%%%%%%%%%%%%%%
\section{Asymptotic estimates of the number of Eulerian circuits}
In what follows we suppose that undirected graph $G$  has no loops and multiple edges, i.e. 
\begin{equation}\label{G_simple}
	G \text{ is a simple graph.}
\end{equation}
We also assume that 
\begin{equation}\label{G_even}
	\text{ all vertices of $G$ have even degrees. }
\end{equation} 

Define the $n \times n$ matrix $Q$ by 
\begin{equation}
Q_{jk} = 
\left\{
\begin{array}{cl}
-1, & (v_j,v_k)\in EG,\\
\phantom{-}d_j,& j = k,\\
\phantom{-}0,  & \text{ otherwise }
\end{array}\right.,
\end{equation}
where $n = |VG|$ and $d_j$ is the degree of the vertex $v_j \in VG$. 
The matrix $Q = Q(G)$ is called the Laplacian matrix of the graph $G$.
The eigenvalues $\lambda_0 \leq \lambda_1 \leq \ldots \leq \lambda_{n-1}$ of the matrix $Q$ are always non-negative real numbers and $\lambda_0 = 0$. 
The eigenvalue $\lambda_1$ is called the algebraic connectivity of the graph $G$.
(For more information about the spectral properties of the Laplace matrix see, for example, \cite{Fiedler1973} and \cite{Mohar1991}.)

According to the Kirchhoff's Matrix-Tree-Theorem, see \cite{Kirchoff1847}, we have that
\begin{equation}\label{Theorem_KMTT}
	t(G) = \frac{1}{n}\lambda_1\lambda_2\cdots\lambda_{n-1},
\end{equation}
where $t(G)$ denotes the number of spanning trees of the graph $G$.

Let $p \geq 1$ be a real number and $\vec{x}\in \mathbb{R}^n$. We use notation 
\begin{equation}
	\left\|\vec{x}\right\|_p = \left(\sum\limits_{j=1}^{n} |x_j|^p \right)^{1/p}.
\end{equation}
For $p = \infty$ we have the maximum norm
\begin{equation}
	\left\|\vec{x}\right\|_\infty = \max_j|x_j|.
\end{equation}
The matrix norm corresponding to the $p$-norm for vectors is
\begin{equation}
	\left\|A\right\|_p = \sup_{\vec{x}\neq 0} \frac{\left\|A\vec{x}\right\|_p}{\left\|\vec{x}\right\|_p}.
\end{equation}
We denote by $\left\|A\right\|_{HS}$ the Hilbert–Schmidt norm of the matrix $A$.
\begin{equation}
 \left\|A\right\|_{HS} = \sqrt{\sum\limits_{j=1}\limits^{n}\sum\limits_{k=1}\limits^{n} |A_{jk}|^2}
\end{equation}

If $f$ is bounded both above and below by $g$ asymptotically, we use the notation
\begin{equation}\label{Big_Theta}
	f(n) = \Theta_{k_1,k_2}\left(g(n)\right),
\end{equation}
which implies as $n\rightarrow \infty$, eventually
\begin{equation}\label{condition_Big_Theta}
	k_1|g(n)| \leq |f(n)| \leq k_2|g(n)|.
\end{equation} 
When functions $f$ and $g$ depend not only on $n$, but also on other parameters $\vec{\xi}$, we use notation (\ref{Big_Theta}) 
meaning that condition (\ref{condition_Big_Theta}) holds uniformly for all possible values of $\vec{\xi}$.

	The main result of the present work is the following theorem.
\begin{Theorem}\label{main}
	Let matrix $Q$ be the Laplacian matrix of graph $G$ with $n$ vertices.	Let conditions (\ref{G_simple}), (\ref{G_even}) hold and 
	%\begin{equation}\label{condition_main}
	%	\left\|(Q+J)^{-1}\right\|_{\infty} \leq \frac{1}{\sigma n},
	%\end{equation}  
	%where $J$ denotes the matrix with every entry $1$. 
	the algebraic connectivity	$\lambda_1 \geq \sigma n$ for some $\sigma > 0$. 
%	for some $\sigma > 1/2$ the degree of each vertex of $G$ at least $\sigma n$.
	Then as $n\rightarrow \infty$
	\begin{equation}\label{main_eq}
		Eul(G) =\Theta_{k_1,k_2} 
		\left(
		  2^{E-\frac{n-1}{2}} \pi^{-\frac{n-1}{2}}\sqrt{t(G)} \,
		 \prod\limits_{j=1}^n \left( \frac{d_j}{2} -1 \right)!
		\right),
	\end{equation}
	where $E = |EG|$, $d_j$ is the degree of the vertex $v_j$,  $t(G)$ denotes the number of spanning trees of the graph $G$ 
	and constants $k_1,k_2 > 0$ depend only on $\sigma$.   
\end{Theorem}

\noindent
{\bf Remark 2.1.} We can replace condition $\lambda_1 \geq \sigma n$ for some $\sigma > 0$ in
Theorem \ref{main} by the condition that for some $\sigma > 1/2$ the degree of each vertex of $G$ at least $\sigma n$.

For the complete graph $K_n$ one can show that $\lambda_1 = n$ and $t(K_n) = n^{n-2}$. 
%The following theorem is the variation of \cite{Brendan1995}.
\begin{Theorem}{(variation of Theorem 4 of \cite{Brendan1995})}\label{TheoremMcc}
	As $n\rightarrow \infty$ with $n$ odd
	\begin{equation}\label{Brendanf}
		Eul(K_n) =  2^{\frac{(n-1)^2}{2}} \pi^{-\frac{n-1}{2}} n^{\frac{n-2}{2}}
		  \left(\left( \frac{n-1}{2} -1 \right)!\right)^n \Big(1 + O(n^{-1/2+\varepsilon})\Big)
	\end{equation}
	for any $\varepsilon>0$.
\end{Theorem}

In fact, Theorem \ref{TheoremMcc} is stronger than Theorem \ref{main} in the case of $G = K_n$. 
However, the asymptotic estimate of Theorem \ref{main} holds for considerably broader class of graphs.
  
%%%%%%%%%%%%%%%%%%%%%%%%%%%%%%%%%%%%%%%%%%%%%%%%%%%%%%%%%%%%%%%%%%%%%%%%%%%%%%%%%%%%%%%%%%%%%%%%%%%%%%%%%%%%%%%%%%%%%%%%%
%%%%%%%%%%%%%%%%%%%%%%%%%%%%%%%%%%%%%%%%%%%%%%%%%%%%%%%%%%%%%%%%%%%%%%%%%%%%%%%%%%%%%%%%%%%%%%%%%%%%%%%%%%%%%%%%%%%%%%%%%
%%%%%%%%%%%%%%%%%%%%%%%%%%%%%%%%%%%%%%%%%%%%%%%%%%%%%%%%%%%%%%%%%%%%%%%%%%%%%%%%%%%%%%%%%%%%%%%%%%%%%%%%%%%%%%%%%%%%%%%%%
%%%%%%%%%%%%%%%%%%%%%%%%%%%%%%%%%%%%%%%%%%%%%%%%%%%%%%%%%%%%%%%%%%%%%%%%%%%%%%%%%%%%%%%%%%%%%%%%%%%%%%%%%%%%%%%%%%%%%%%%%

\section{Some basic properties of the Laplacian matrix}
Consider the graph $G$ such that conditions (\ref{G_simple}), (\ref{G_even}) hold. 
The Laplacian matrix $Q$ has the eigenvector $[1,1,\ldots,1]^T$, corresponding to the eigenvalue $\lambda_0 = 0$.
%For our purpose it is convenient to denote by $L$ the orthogonal complement to $[1,1,\ldots,1]^T$. 
We use notation $\hat{Q} = Q + J$, where $J$ denotes the matrix with every entry $1$. Note that $Q$ and $\hat{Q}$ 
have the same set of eigenvectors and eigenvalues, except for the eigenvalue corresponding to the eigenvector $[1,1,\ldots,1]^T$, which equals
$0$ for $Q$ and $n$ for $\hat{Q}$.

Since the spectral norm  is bounded above by any matrix norm we get that 
\begin{equation}\label{lambda<n}
	\lambda_{n-1} = ||Q||_2 \leq ||\hat{Q}||_2 \leq ||\hat{Q}||_1 = \max_{j}{\sum\limits_{k=1}^{n}} |\hat{Q}_{jk}| = n.
\end{equation}

We denote by $G_r$ the graph which arises from $G$ by removing vertices $v_1, v_2, \ldots, v_r$ and all adjacent edges.
%The following lemma proved in  \cite{Fiedler1973}. 
\begin{Lemma}\label{Lemma_connectivity}
	Let condition (\ref{G_simple}) holds for graph $G$  with $n$ vertices. Then 
	\begin{equation}\label{lambda<d}
		\lambda_1(G) \leq \frac{n}{n-1} \min_j d_j,
	\end{equation}
	\begin{equation}\label{lambda>d}
		\lambda_1(G) \geq 2 \min_j d_j - n + 2,
	\end{equation}
	\begin{equation}\label{lambda_Gr}
		\lambda_1(G_r) \geq \lambda_1(G) - r,
	\end{equation}
	where $\lambda_1(G)$ is the algebraic connectivity of $G$ and $d_j$ is the degree of the vertex $v_j \in VG$. 
\end{Lemma}
The proof of Lemma \ref{Lemma_connectivity} can be found in \cite{Fiedler1973}.

\begin{Lemma}\label{Lemma_normQ}
	Let condition (\ref{G_simple}) hold and 
	%\begin{equation}\label{condition_main}
	%	\left\|(Q+J)^{-1}\right\|_{\infty} \leq \frac{1}{\sigma n},
	%\end{equation}  
	%where $J$ denotes the matrix with every entry $1$. 
	the algebraic connectivity	$\lambda_1 \geq \sigma n$ for some $\sigma > 0$. 
	Then there is a constant $c_{\infty}>0$ depending only on $\sigma$ such that 
	\begin{equation}\label{normQ}
		||\hat{Q}^{-1}||_1 = ||\hat{Q}^{-1}||_\infty \leq \frac{c_{\infty}}{n}. 
	\end{equation}
\end{Lemma}
\begin{Proof} {\it Lemma \ref{Lemma_normQ}.} 
	We consider $\vec{x}\in \mathbb{R}^n$ such that $||\vec{x}||_\infty = 1$. 
	For simplicity, we assume that $|x_1| = 1$.
	We denote by $J_{\sigma}$ 
	the set of the indices $j$ such that $|x_j|\geq \sigma/8$. 
	%Since the algebraic connectivity 	$\lambda_1 \geq \sigma n$, using (\ref{lambda<d}), 
	%we get that 
	%\begin{equation}\label{Lemma_normQ_d1}
%		d_1 \geq \min_j d_j \geq \frac{n-1}{n} \sigma n \geq \sigma n/2.
%	\end{equation}

	In the case of $|J_{\sigma}| \geq \sigma n/4$ we have that
	\begin{equation}
		||\vec{x}||_2 \geq \sqrt{\frac{\sigma^2}{64}\sigma n/4}.
	\end{equation}
	Since the algebraic connectivity 	$\lambda_1 \geq \sigma n$, we get that
	\begin{equation}\label{Lemma_normQ_J<sigma}
		\sqrt{n ||\hat{Q}\vec{x}||^2_\infty} 
		\geq
		||\hat{Q}\vec{x}||_2 
		\geq
		\lambda_1||\vec{x}||_2
		\geq
		\sigma n||\vec{x}||_2
		\geq \sigma n \sqrt{\frac{\sigma^3 n}{256}}
	\end{equation}
	In the case of $|J_{\sigma}| \leq \sigma n/4$  we have that
	\begin{equation}\label{Lemma_normQ_J>sigma}
		\begin{aligned}
			||\hat{Q}\vec{x}||_\infty
			&\geq (d_1+1)|x_1| - \sum\limits_{j=2}\limits^{n}|x_j| \geq \\
			&\geq d_1+1 - \sum\limits_{j\in J_\sigma}|x_j| - \sum\limits_{j \notin J_\sigma}|x_j| \geq \\
			&\geq d_1+1 - \sigma n /{4} - n \sigma/8.
		\end{aligned}
	\end{equation}
	Using again $\lambda_1 \geq \sigma n$ and (\ref{lambda<d}) we get that 
	\begin{equation}\label{Lemma_normQ_d1}
		d_1 \geq \min_j d_j \geq \frac{n-1}{n} \sigma n \geq \sigma n/2.
	\end{equation}
	Combining (\ref{Lemma_normQ_J<sigma}), (\ref{Lemma_normQ_J>sigma}) and (\ref{Lemma_normQ_d1}) we obtain that
	\begin{equation}
	||\hat{Q}\vec{x}||_\infty \geq c_{\infty}^{-1} n ||\vec{x}||_\infty,
	\end{equation}
	for some constant $c_{\infty}>0$ depending only on $\sigma.$  
\end{Proof}

The following lemmas will be applied to estimate the
determinant of a matrix close to the identity matrix $I$.
\begin{Lemma}\label{Lemma_matrix1}
	Let $X$ be an $n\times n$ matrix  such that $\left\|X\right\|_2 < 1$. Then for fixed $m \geq 2$
	\begin{equation}
		\det(I+X) = \exp \left(  \sum\limits_{r=1}^{m-1} \frac{(-1)^{r+1}}{r}\, tr (X^r) + E_m(X)   \right),
	\end{equation}
 	where \text{tr} is the trace function and
	\begin{equation}
		|E_m(X)|\leq \frac{n}{m}\,\frac{\left\|X\right\|_2^m}{1-\left\|X\right\|_2}.
	\end{equation}
\end{Lemma}
Lemma \ref{Lemma_matrix1} was also formulated and proved in \cite{Brendan1995}.

\begin{Lemma}\label{Lemma_matrix2}
	Let the assumptions of Lemma \ref{Lemma_matrix1} hold and all eigenvalues of $X$ are non-negative real numbers. Then
\begin{equation}
\det(I-X) \geq \exp \left( -\frac{tr(X)}{1-\left\|X\right\|_2} \right).
\end{equation}
\end{Lemma}

\begin{Proof} {\it Lemma \ref{Lemma_matrix2}.} 
	Using Lemma \ref{Lemma_matrix1} we get that
	\begin{equation}
		\det(I-X) = \exp \left(  \sum\limits_{r=1}^{\infty} \frac{(-1)^{r+1}}{r}\, (-1)^r tr(X^r)    \right)
	\end{equation}
	Since all eigenvalues of $X$ are non-negative real numbers  
	\begin{equation}
		0\leq tr (X^r) \leq tr (X) \left\|X\right\|_2^{r-1}.
	\end{equation}
	Hence 
	\begin{equation}
		\det(I-X) 
		\geq 
		\exp \left(  - \sum\limits_{r=1}^{\infty} \frac{tr (X)}{r} \left\|X\right\|_2^{r-1}    \right)
		\geq 
		\exp \left( -\frac{tr(X)}{1-\left\|X\right\|_2} \right).
	\end{equation}
\end{Proof}

\begin{Lemma}\label{Lemma_minor}
		Let the assumptions of Lemma \ref{Lemma_normQ} hold. Then there is a constant $c_1>0$ depending only on $\sigma$ 
	such that 
	\begin{equation}\label{eq_Lemma_minor}
		|\det M_{11}| \leq c_1 \frac{\det \hat{Q}}{n},	
	\end{equation} 
	where $M_{11}$ denotes the $(n - 1)\times(n - 1)$ matrix that results from deleting the first row and the first column of 
	$\hat{Q} = Q + J$.
\end{Lemma}
\begin{Proof} {\it Lemma \ref{Lemma_minor}.} 
 Since the algebraic connectivity 	$\lambda_1 \geq \sigma n$, using (\ref{lambda<d}), 
	we get the following estimate for the degree $d_k$ of the vertex $v_k \in VG.$ 
	\begin{equation}\label{Lemma_minor_d1}
		d_k \geq \min_j d_j \geq \frac{n-1}{n} \sigma n \geq \sigma n/2.
	\end{equation}
	Consider the $n\times n$ matrix $X$ such that
	\begin{equation}
		X_{jk} = 
		\left\{
			\begin{array}{cl}
			\frac{1}{d_1+1},& \text{ if  }  (v_1, v_j) \notin EG, (v_1, v_k) \notin EG \text{ and } j,k\neq 1,\\
			0 ,& \text{ otherwise.} 
			\end{array}
		\right.
	\end{equation}
	After performing one step of the Gaussian elimination for $\hat{Q}+X$, we obtain that
	\begin{equation}\label{Lemma_minor_XM}
		\det(\hat{Q} + X) = (d_1+1) \det M_{11},
	\end{equation} 
	Since the spectral norm is bounded above by any matrix norm, 
	%using (\ref{Lemma_minor_d1}), 
	we get that
	\begin{equation}\label{Lemma_minor_X2}
		||X||_2 \leq ||X||_1 \leq \frac{n}{d_1+1} \leq \frac{2}{\sigma}.
	\end{equation}
	Since $\lambda_1 \geq \sigma n$, taking  into account (\ref{lambda<n}), we obtain that
	\begin{equation}\label{Lemma_minor_XQ2}
		||X\hat{Q}^{-1}||_2 \leq ||X||_2 ||\hat{Q}^{-1}||_2 \leq \frac{2}{ \sigma \lambda_1} \leq \frac{2}{ \sigma^2 n} 
	\end{equation}
	Combining Lemma \ref{Lemma_matrix1} with (\ref{Lemma_minor_XQ2}), we get that as $n\rightarrow \infty$
	\begin{equation}\label{Lemma_minor_I+X}
		\det\left({I+X\hat{Q}^{-1}}\right) = \exp\left( tr \left( X\hat{Q}^{-1}\right) + E_2\left(X\hat{Q}^{-1}\right)\right) 
		\leq \exp\left( n\frac{2}{ \sigma^2 n} + O(n^{-1})\right).
	\end{equation}
	From (\ref{Lemma_minor_XM}) and  (\ref{Lemma_minor_I+X}) we have that as $n\rightarrow \infty$
	\begin{equation}\label{Lemma_minor_last}
			(d_1+1) \det M_{11} =
		%	\det\left(\hat{Q} + X\right) =
			 \det\left({I+X\hat{Q}^{-1}}\right)\det\hat{Q} 
			 \leq
			 \det\hat{Q} \exp\left( 2/\sigma^2 + O(n^{-1})\right).
	\end{equation}
	Since $\hat{Q}$ is positive definite,
	using (\ref{Lemma_minor_d1}) in (\ref{Lemma_minor_last}), we obtain (\ref{eq_Lemma_minor}). 
\end{Proof}

\begin{Lemma}\label{Lemma_Gr}
	Let the assumptions of Lemma \ref{Lemma_normQ} hold. 
	Let $G_r$ be the graph which arises from $G$ by removing vertices $v_1, v_2, \ldots, v_r$ and all adjacent edges. 
	 Then there is a constant $c_2>0$ depending only on $\sigma$ such that for any $\varepsilon \in (0,1)$ and $r\leq n^{\varepsilon}$
	 \begin{equation}\label{eq_Lemma_Gr}
	 	\det{\hat{Q}(G_r)} \geq \frac{\det\hat{Q}(G)}{\left(c_2 n\right)^r}.
	 \end{equation}
\end{Lemma}
\begin{Proof} {\it Lemma \ref{Lemma_Gr}.}  
	We give first a proof for the case of $r=1$. For our purpose it is convenient to use notations $\hat{Q} = \hat{Q}(G)$ and $\hat{Q}_1 = \hat{Q}(G_1)$.
	Note that the matrix $M_{11}$ that results from deleting the first row and the first column of 
	$\hat{Q}$ coincides with the matrix $\hat{Q}_1$ with the exception of the diagonal elements. 
	In a similar way as (\ref{Lemma_minor_XM}) we get that
	\begin{equation}\label{Lemma_Gr_G1}
		\det(\hat{Q} + \Omega+ X ) = (d_1+1) \det \hat{Q}_1,
	\end{equation} 
	where $X$ is such that  
	\begin{equation}
		X_{jk} = 
		\left\{
			\begin{array}{cl}
			\frac{1}{d_1+1},& \text{ if } (v_1, v_j) \notin EG, (v_1, v_k) \notin EG \text{ and } j,k\neq 1,\\
			0 ,& \text{ otherwise.} 
			\end{array}
		\right.
	\end{equation}
	 and $\Omega$ is such diagonal matrix that  
	\begin{equation}
		\Omega_{jj} = 
		\left\{
			\begin{array}{cl}
			1,& \text{ if } (v_1, v_j) \in EG,\\
			0 ,& \text{ otherwise.} 
			\end{array}
		\right.
	\end{equation}
		Taking  into account (\ref{Lemma_minor_X2}), we have that
	\begin{equation}
		||\Omega+ X||_2 \leq ||\Omega||_2 + ||X||_2 \leq 1+\frac{2}{\sigma} \leq \frac{3}{\sigma}.
	\end{equation}
	In a similar way as (\ref{Lemma_minor_XQ2}) we get that
	\begin{equation}\label{Lemma_Gr_XQ2}
		||(\Omega+X)\hat{Q}^{-1}||_2 \leq ||\Omega+X||_2 ||\hat{Q}^{-1}||_2 
		\leq \frac{3}{\sigma\lambda_1} \leq \frac{3}{ \sigma^2 n}
	\end{equation}
	Combining Lemma \ref{Lemma_matrix1} with (\ref{Lemma_Gr_XQ2}), we get that as $n\rightarrow \infty$
	\begin{equation}\label{Lemma_Gr_I+X}
		\begin{array}{c}
		\displaystyle
		\det\left({I+(\Omega+X)\hat{Q}^{-1}}\right) 
		= 
		\exp\left( tr \left((\Omega+X)\hat{Q}^{-1}\right) + E_2\left((\Omega+X)\hat{Q}^{-1}\right)\right) 
		\geq
		\\\displaystyle
		\geq 
		\exp\left( -n\frac{3}{ \sigma^2 n} + O(n^{-1})\right).
		\end{array}
	\end{equation}
	From (\ref{Lemma_Gr_G1}) and (\ref{Lemma_Gr_I+X}) we have that as $n \rightarrow \infty$
		\begin{equation}\label{Lemma_Gr_last}
			(d_1+1) \det \hat{Q}_{1} =
		%	\det\left(\hat{Q} + X\right) =
			 \det\left({I+(\Omega+X)\hat{Q}^{-1}}\right)\det\hat{Q} 
			 \geq
			 \det\hat{Q} \exp\left( -3/\sigma^2 + O(n^{-1})\right).
	\end{equation}
	Since $d_1+1 \leq n$ we get (\ref{eq_Lemma_Gr}) for the case of $r = 1$.
	
	Taking  into account  (\ref{lambda_Gr}) and using $r$ times  (\ref{Lemma_Gr_last}) we get (\ref{eq_Lemma_Gr})
	for the general case.
\end{Proof}

According to (\ref{Theorem_KMTT}), we have that 
 \begin{equation}\label{t(G)hat}
 	t(G) = \frac{1}{n}\lambda_1\lambda_2\cdots\lambda_{n-1} = \frac{\det{\hat{Q}}}{n^2},
 \end{equation}
where $t(G)$ denotes the number of spanning trees of the graph $G$.
\begin{Lemma}\label{Lemma_Trees>d}
	Let the assumptions of Lemma \ref{Lemma_normQ} hold.  Then for some $c_3> 0$ depending only on $\sigma$ the number of spanning trees of the graph $G$ 
	with maximum degree greater than $d$ is less then $c_3^n\, \det{\hat{Q}}/d!$ for all $d \geq 0$.
\end{Lemma}
\begin{Proof} {\it Lemma \ref{Lemma_Trees>d}.} 
	According to Lemma 5 of \cite{Brendan1995} the number of labelled trees on $n$ vertices with first vertex having degree greater than
$d$ is less than $2n^{n-2}/d!$ for all $d \geq 0$. We have that
	\begin{equation}
		\det{\hat{Q}} \geq \lambda_1^n \geq (\sigma n)^n.
	\end{equation}
To complete proof it remains to note that the number of spanning trees with maximum degree greater than 
$d$ in $G$ does not exceed the number of such spanning trees in the complete graph with $n$ vertices.
\end{Proof}

Consider a spanning tree $T$ of the graph $G$. 
We denote by $G_T$ the graph which arises from $G$ by removing all edges of the tree $T$.
\begin{Lemma}\label{Lemma_G_T}
			Let the assumptions of Lemma \ref{Lemma_normQ} hold. Let $T$ be the spanning tree of $G$ 
		with the maximum degree at most $\sigma n /4$. Then the algebraic connectivity $\lambda_1(G_T) \geq \sigma n /2$ 
		and 
		\begin{equation}
			\det{\hat{Q}(G_T)} \geq c_4 \det{\hat{Q}(G)} 
		\end{equation}
		for some $c_4>0$ depending only on $\sigma$.
\end{Lemma}
\begin{Proof} {\it Lemma \ref{Lemma_G_T}.}
	Note that 
	\begin{equation}
		\hat{Q}(G_T) = \hat{Q}(G) - Q(T).
	\end{equation}
	Since the spectral norm  is bounded above by any matrix norm 
	and the maximum degree of vertex of $T$ at most $\sigma n /4$
	we get that
	\begin{equation}\label{Lemma_G_T_Q_T}
		||Q(T)||_2 \leq ||Q(T)||_1 \leq \sigma n /2.
	\end{equation}
	Therefore, since the algebraic connectivity $\lambda_1(G) \geq \sigma n$, we have that
	\begin{equation} 
		\lambda_1(G_T) \geq \lambda_1(G) - ||Q(T)||_2 \geq \sigma n /2
	\end{equation}
	and 
	\begin{equation}\label{Lemma_G_T_I-X}
		\det \hat{Q}(G_T) = \det \hat{Q}(G) \det (I - X),
	\end{equation}
	where $X = Q(T)\hat{Q}(G)^{-1}$. 
	Note that $\hat{Q}(G)$ is the matrix of positive definite quadratic form and 
	 $Q(T)$ is the matrix of quadratic form with non-negative eigenvalues.
	 Considering the basis in which both matrices are diagonal, we have that
	\begin{equation}\label{Lemma_G_T_trQQ}
		tr(Q(T)\hat{Q}(G)^{-1}) \leq tr(Q(T)) ||\hat{Q}(G)^{-1}||_2
	\end{equation}
	and
	\begin{equation}
		\text{ all eigenvalues of } Q(T)\hat{Q}(G)^{-1} \text{ are non-negative. }
	\end{equation}
	Using again the fact that the algebraic connectivity $\lambda_1(G) \geq \sigma n$ and (\ref{Lemma_G_T_Q_T}) we get that
	\begin{equation}\label{Lemma_G_TXX}
		||X||_2 \leq ||Q(T)||_2 ||\hat{Q}(G)^{-1}||_2 \leq \frac{\sigma n}{2} \frac{1} {\sigma n} = \frac{1}{2}.
	\end{equation}	
	Since $T$ is the spanning tree 	$tr(Q(T)) = 2 (n - 1)$. Using (\ref{Lemma_G_T_trQQ}), we get that
	\begin{equation}\label{Lemma_G_TtrX}
		tr(X) \leq tr(Q(T)) ||\hat{Q}(G)^{-1}||_2 \leq 2 n \frac{1} {\sigma n} = \frac{2}{\sigma}.
	\end{equation}
	To complete the proof it remains to combine Lemma \ref{Lemma_matrix2} with (\ref{Lemma_G_T_I-X}), (\ref{Lemma_G_TXX}) and (\ref{Lemma_G_TtrX}).
\end{Proof} 

\begin{Lemma}\label{Lemma_struct}
	Let $a>0$ and the assumptions of Lemma \ref{Lemma_normQ} hold. Then for any set  $A \subset VG$ such, that
	  $|A| \geq an$, there is a function $h: VG \rightarrow \mathbb{N}_0$, 
	  having following properties: 
	   \begin{equation}
	  	h(v) = 0, \text { if } v \in A, \ \ \	h(v) \leq H,\text { for any }  v \in VG,
	  \end{equation}
	 	\begin{equation}\label{Eq_3_30}
	  	 \Big|\left\{w \in VG \ | \ (w,v)\in EG \text{ and } h(w)<h(v)\right\}\Big| \geq \alpha n, \text { if } v \notin A, 
	  \end{equation}
	  where constants $H, \alpha > 0$ depend only on $a$ and  $\sigma$.
\end{Lemma}
\begin{Proof} {\it Lemma \ref{Lemma_struct}.}
	At first, we construct the set $A_1 = \left\{v \in VG \ | \ h(v) = 1 \right\}$, having property (\ref{Eq_3_30}). 
	
	 If $|A| > n - \sigma n/4$, then let  $A_1 = \left\{v \in VG \ | \ v \notin A\right\}.$
	 Taking into account (\ref{Lemma_minor_d1}), we get that property (\ref{Eq_3_30}) hold for $\alpha = \sigma/4$. In this case $H = 1$. 
	 
	 For $|A| \leq n - \sigma n/4$ define $\vec{x} \in \mathbb{R}^n$ such that 
	 \begin{equation}
	 	x_j = 
	 	\left\{
			\begin{array}{cl}
			1 - |A|/n,& \text{  }  v_j \in A,\\
			-|A|/n ,& \text{ }  v_j \notin A. 
			\end{array}
		\right.
	 \end{equation}
	 Since $x_1+x_2+\ldots+x_n =0$ 
	 \begin{equation}\label{Eq_3_32}
	 		\vec{x}^T Q \vec{x} = \vec{x}^T \hat{Q} \vec{x} \geq \lambda_1 \|\vec{x}\|_2^2 \geq \lambda_1 |A|\, \left(\frac{n - |A|}{n}\right)^2 
	 		\geq \sigma n \, a n \left(\sigma/4\right)^2 = \frac{a\sigma^3 n^2}{16}.
	 \end{equation}
	 On the other hand, 
	 \begin{equation}
	 	\vec{x}^T Q \vec{x} = \sum\limits_{(v_j,v_k)\in EG} (x_j - x_k)^2,
	 \end{equation} 
	 which is equal to the number of edges $(v,w)\in EG$, where $v\in A, w \notin A.$ 
	 We denote $A_1$ the set of vertices $w \notin A$, having at least $\alpha n$ adjacent vertices in $A$, 
	 where $\alpha = \frac{1}{32}a\sigma^3$. 
	 \begin{equation}\label{Eq_3_34}
	 			\vec{x}^T Q \vec{x} \leq n |A_1| + \alpha n |VG|.
	 \end{equation}
	 Combining (\ref{Eq_3_32}) and (\ref{Eq_3_34}), we get that $|A_1| \geq \alpha n.$
	 
	 We make further construction of the function $h$ inductively, using for the $k$-th step the set $A^{(k)} = A\cup A_1\cup \ldots \cup A_k.$ 
	 The number of steps does not exceed $1/\alpha$ as $|A_k| \geq \alpha n$ for each step, perhaps with the exception of the last one.
\end{Proof}

%%%%%%%%%%%%%%%%%%%%%%%%%%%%%%%%%%%%%%%%%%%%%%%%%%%%%%%%%%%%%%%%%%%%%%%%%%%%%%%%%%%%%%%%%%%%%%%%%%%%%%%%%%%%%%%%%%%%%%%%%
%%%%%%%%%%%%%%%%%%%%%%%%%%%%%%%%%%%%%%%%%%%%%%%%%%%%%%%%%%%%%%%%%%%%%%%%%%%%%%%%%%%%%%%%%%%%%%%%%%%%%%%%%%%%%%%%%%%%%%%%%
%%%%%%%%%%%%%%%%%%%%%%%%%%%%%%%%%%%%%%%%%%%%%%%%%%%%%%%%%%%%%%%%%%%%%%%%%%%%%%%%%%%%%%%%%%%%%%%%%%%%%%%%%%%%%%%%%%%%%%%%%
%%%%%%%%%%%%%%%%%%%%%%%%%%%%%%%%%%%%%%%%%%%%%%%%%%%%%%%%%%%%%%%%%%%%%%%%%%%%%%%%%%%%%%%%%%%%%%%%%%%%%%%%%%%%%%%%%%%%%%%%%

\section{The result expressed as an integral}
The reasoning below is similar to the arguments of Section 2 of \cite{Brendan1995}

An Eulerian orientation of $G$ is an orientation of its edges with the property that for every
vertex both the in-degree and the out-degree are equal. Any Eulerian circuit
induces an Eulerian orientation by orienting each edge in accordance with its direction of
traversal.

A directed tree with root $v$ is a connected directed graph $T$ such that $v \in VT$ has out-
degree zero, and each other vertex has out-degree one. Thus, $T$ is a tree which has each edge
oriented towards $v$.

Let $D$ be a directed graph with $n$ vertices, and let $v \in VD$. A directed spanning tree of $D$
with root $v$ is a spanning subgraph of $D$ which is a directed tree with root $v$.

The following famous theorem, sometimes called the BEST Theorem, is due to de Bruijn,
van Aardenne-Ehrenfest, Smith, and Tutte \cite{Aardenne1951, Smith1941}.
\begin{Theorem}\label{BEST}
	Let $D$ be a directed graph with vertices $v_1, v_2, \ldots , v_n$. Suppose that there are
numbers $d_1, d_2, \ldots , d_n$ such that, for every vertex $v_r$, both the in-degree and the out-degree of
$v_r$ are equal to $d_r$. Let $t_r = t_r (D)$ be the number of directed spanning trees of $D$ rooted at $v_r$.
Then $t_r$ is independent of $r$, and
	\begin{equation}
		Eul(D) = t_r\prod\limits_{j=1}\limits^{n}(d_j-1)!.
	\end{equation}
\end{Theorem}

	Consider the undirected graph $G$ with $n$ vertices such that conditions (\ref{G_simple}), (\ref{G_even}) hold.
	Note that for every spanning tree $T$ of the graph $G$ and any vertex $v_r\in VG$ there is only one orientation of the edges of $T$
 such that we obtain a directed tree with root $v_r$. We denote by ${\cal T}_r$ the set of directed trees with root $v_r$ obtained in such a way.
 For $T\in\bigcup\limits_{r=1}\limits^n {\cal T}_r$ denote by $EO(T)$  the number of Eulerian orientations of $G$ that 
 the corresponding graphs contain $T$. 

From Theorem \ref{BEST}  in the case of a graph $D$ corresponding to Eulerian orientation of the graph $G$ we find that
\begin{equation}
	Eul(D) = t_r(D)\prod\limits_{j=1}\limits^{n}\left(\frac{d_j}{2}-1\right)!,
\end{equation}
 where $d_j$ is the degree of the vertex $v_j \in VG$. Let denote by ${\cal EO}$ the set of all graphs corresponding to Eulerian orientations 
 of the graph $G$. Grouping Eulerian circuits according to the induced orientations, we obtain that
\begin{equation}
	Eul(G) = \sum\limits_{D \in {\cal EO}} Eul(D) = \prod\limits_{j=1}\limits^{n}\left(\frac{d_j}{2}-1\right)! \sum\limits_{D \in {\cal EO}} t_r(D)
\end{equation}
for any fixed natural number $r\leq n$.

Regrouping the terms of the final summation according to the directed subtrees rooted at $v_r$,
we find that
\begin{equation}
	Eul(G) = \prod\limits_{j=1}\limits^{n}\left(\frac{d_j}{2}-1\right)! \sum\limits_{T \in {\cal T}_r} EO(T).
\end{equation}

For $n \geq 1$ and $R \geq 0$  we use notation  $U_n(R) = \{(x_1, x_2, \ldots, x_n)\ | \ |x_i| < R \text{ for all } i\}$. The value of
$EO(T)$ is the constant term in
\begin{equation}
	\prod\limits_
		{(v_j,v_k) \in EG}
		({x_j}^{-1} x_k + {x_k}^{-1} x_j) 
	\prod\limits_{(v_j,v_k)\in ET} \frac{x_k^{-1} x_j}{({x_j}^{-1} x_k + {x_k}^{-1} x_j)},
\end{equation}
which we can extract via Cauchy’s Theorem using the unit circle as a contour for each variable.
Making the substitution $x_j = e^{i\theta_j}$ for each $j$, we find that
\begin{equation}\label{EulS}
Eul(G) = \prod\limits_{j=1}^n \big( \frac{d_j}{2} -1 \big)! \, 2^{|EG|-n+1} \pi^{-n} S,
\end{equation}
where
\begin{equation}\label{S}
	S = \int\limits_{U_n(\pi/2)} 
 % \left( 
     \prod\limits_{(v_j,v_k)\in EG} \cos \Delta_{jk} \sum \limits_{T \in {\cal T}_r} 
     \prod \limits_{(v_j,v_k)\in ET} (1+i\tan \Delta_{jk})  
 % \right) 
 \
 d\vec{\theta},
\end{equation}
having put $\Delta_{jk} = \theta_j - \theta_k$ and using the fact that the integrand is unchanged by the substitutions
 $\theta_j \rightarrow \theta_j + \pi$ if condition (\ref{G_even}) holds.

We approach the integral by first estimating it in the region that would turn out to be the
asymptotically significant one. Then we bound the integral over the remaining regions and
show that it is vanishingly small in comparison with the significant part.

%%%%%%%%%%%%%%%%%%%%%%%%%%%%%%%%%%%%%%%%%%%%%%%%%%%%%%%%%%%%%%%%%%%%%%%%%%%%%%%%%%%%%%%%%%%%%%%%%%%%%%%%%%%%%%%%%%%%%%%%%
%%%%%%%%%%%%%%%%%%%%%%%%%%%%%%%%%%%%%%%%%%%%%%%%%%%%%%%%%%%%%%%%%%%%%%%%%%%%%%%%%%%%%%%%%%%%%%%%%%%%%%%%%%%%%%%%%%%%%%%%%
%%%%%%%%%%%%%%%%%%%%%%%%%%%%%%%%%%%%%%%%%%%%%%%%%%%%%%%%%%%%%%%%%%%%%%%%%%%%%%%%%%%%%%%%%%%%%%%%%%%%%%%%%%%%%%%%%%%%%%%%%
%%%%%%%%%%%%%%%%%%%%%%%%%%%%%%%%%%%%%%%%%%%%%%%%%%%%%%%%%%%%%%%%%%%%%%%%%%%%%%%%%%%%%%%%%%%%%%%%%%%%%%%%%%%%%%%%%%%%%%%%%

\section{ The dominant part of the integral }
%For our purpose it is convenient to denote by $L$ the orthogonal complement to $[1,1,\ldots,1]^T$. 
In what follows, we fix some small constant  $\varepsilon>0$. Define
\begin{equation}
	V_0 = \{\vec{\theta}\in U_n(\pi/2) :\ |\theta_j - \bar{\theta}| \, (\, mod \, \pi) \leq n^{-1/2+\varepsilon}\text{, where } \bar{\theta}=\frac{\theta_1+\ldots+\theta_n}{n}\},
\end{equation}
and let $S_0$ denote the contribution  to $S$ of $\vec{\theta}\in V_0$. 
%\begin{equation}\label{S_0}
%	S_0 = \int\limits_{V_0}
 % \left( 
%     \prod\limits_{(v_j,v_k)\in EG} \cos \Delta_{jk} \sum \limits_{T \in {\cal T}_r} 
 %    \prod \limits_{(v_j,v_k)\in ET} (1+i\tan \Delta_{jk})  
 % \right) 
 \
% d\vec{\theta}.
%\end{equation}
Since the integrand is invariant under uniform translation of all
the $\theta_j$'s mod $\pi$, we can fix $\bar{\theta} = 0$ and multiply it by the ratio 
of its range $\pi$ to the length $n^{-1/2}$ of the vector $\frac{1}{n}[1,1,\ldots,1]^T$. 
Thus we have that
\begin{equation}\label{S_0} 
S_0 = \pi n^{1/2} \int\limits_{L \cap V_0}
%\left( 
     \prod\limits_{(v_j,v_k)\in EG} \cos \Delta_{jk} \sum \limits_{T \in {\cal T}_r} 
     \prod \limits_{(v_j,v_k)\in ET} (1+i\tan \Delta_{jk})  
% \right) 
  \ dL,
\end{equation}
where $L$ denotes the orthogonal complement to the vector $[1,1,\ldots,1]^T$.

The sum over ${\cal T}_r$ in the integrand of (\ref{S_0}) can be expressed as a determinant, according to the following theorem of \cite{Tutte1948}.

\begin{Theorem}\label{Tutte}
Let $w_{jk}$ $(1 \leq j, k \leq n$, $j \neq k$) be arbitrary. Define the $n \times n$ matrix $A$ by
\begin{equation}\label{matrix_Tutte}
A_{jk} = 
\left\{
\begin{array}{cc}
-w_{jk},& \text{ if } j\neq k,\\
\sum_{r\neq j} w_{jr},& \text{ if } k=j
\end{array}
\right.,
\end{equation}
the sum being over $1 \leq r \leq n$ with $r\neq j$. For any $r$ with $1 \leq r \leq n$, let 
 $M_r$ denote the
principal minor of $A$ formed by removing row $r$ and column $r$. Then
\begin{equation}
\det M_r = \sum\limits_{T} \prod\limits_{(v_j, v_k)\in ET} w_{jk},
\end{equation}
where the sum is over all directed trees $T$ with $VT = \{v_1, v_2, \ldots, v_n\}$ and root $v_r$.
\end{Theorem}

\begin{Lemma}\label{Lemma_using_Tutte}
	Let the assumptions of Theorem \ref{main} hold. Let $\hat{Q} = Q + J$, where $J$ denotes the matrix with every entry 1.
	Then for $\vec{\theta} \in V_0$ as $n\rightarrow \infty$
	\begin{equation}\label{eq_Lemma_using_Tutte}
		\sum\limits_{r=1}\limits^{n}\sum\limits_{T \in {\cal T}_r} \prod\limits_{(v_j,v_k)\in ET} (1+i\tan \Delta_{jk})
			=
		\frac{e^{i\vec{\theta}^T Q \vec{\alpha} + \frac{1}{2}tr(\Lambda\hat{Q}^{-1})^2}}{n}
		\det \hat{Q} \Big( 1+O(n^{-1/2+3\varepsilon}) \Big),
	\end{equation}
where $\vec{\alpha}$ denotes the vector composed of the diagonal elements of the matrix $\hat{Q}^{-1}$, 
$\Lambda$ denotes the diagonal matrix whose diagonal elements are equal to the components of the vector $Q\vec{\theta}$. 
\end{Lemma}

\begin{Proof} {\it Lemma \ref{Lemma_using_Tutte}.}
	Define the $n\times n$ matrix $B$ by
	\begin{equation}
		B_{jk} = 
		\left\{
		\begin{array}{ll}
			-\tan\Delta_{jk},& \text{ for } (v_j, v_k) \in EG,\\
			\sum\limits_{l: (v_j, v_l) \in EG} \tan \Delta_{jl},& \text{ for } k=j,\\
			0 & \text{ otherwise }.
			\end{array}
		\right.
	\end{equation}
	Using Theorem \ref{Tutte} with the matrix $A = Q + iB$, we get that 
	\begin{equation}\label{Lemma_using_Tutte_000}
			\sum\limits_{r=1}\limits^{n}\sum\limits_{T \in {\cal T}_r} \prod\limits_{(v_j,v_k)\in ET} (1+i\tan \Delta_{jk}) = 
			\sum\limits_{r=1}\limits^{n} M_r,	
	\end{equation} 	
	where $M_r$ denotes the
	principal minor of $A$ formed by removing row $r$ and column $r$.
	Since the vector $[1,1,\ldots,1]^T$ is the common eigenvector of the matrices $Q$ and $B$, corresponding to the eigenvalue $0$, we find that
	\begin{equation}\label{Lemma_using_Tutte_001}
		\sum\limits_{r=1}\limits^{n} M_r = \frac{\det(\hat{Q} + iB)}{n}.
	\end{equation}
	Note that for $\vec{\theta} \in V_0$
	\begin{equation}\label{Lemma_using_Tutte_Delta}
		\Delta_{jk} = (\theta_j - \bar{\theta}) - (\theta_k - \bar{\theta}) \leq 2 n^{-1/2+\varepsilon}.
	\end{equation}
	Since the spectral norm is bounded above by any matrix norm we get that
	\begin{equation}\label{Lemma_using_Tutte_B}
		||B||_2 \leq ||B||_1 = \max_{j}{\sum\limits_{k=1}^{n}} |B_{jk}| = O(n^{1/2+\varepsilon}).
	\end{equation}
	Let $\Phi = B\hat{Q}^{-1}$. Since the algebraic connectivity 	$\lambda_1 \geq \sigma n$, we get that
	\begin{equation}
		||\Phi||_2 \leq ||B||_2 ||\hat{Q}^{-1}||_2 \leq \frac{1}{\lambda_1}||B||_2 = O(n^{-1/2+\varepsilon}).
	\end{equation}
	 Using Lemma (\ref{Lemma_matrix1}) with the matrix $i\Phi$, we find that as $n\rightarrow\infty$ 
	 \begin{equation}\label{Lemma_using_Tutte_I+Phi}
	 	\det(I+i \Phi) = \exp \left( tr (i \Phi) + \frac{tr (\Phi^2)}{2}  + O(n^{-1/2+3\varepsilon}) \right).
	 \end{equation}
	 Let 
	 \begin{equation}
	 	B = B_{skew} + B_{diag},
	 \end{equation}
	 where $B_{skew}$ is the skew-symmetric matrix and $B_{diag}$ is the diagonal matrix.
	 Since $\hat{Q}$ is the symmetric matrix 
 	 \begin{equation}\label{Bskewsmall}
 	  	tr(B_{skew}\hat{Q}^{-1}) = 0.
 	 \end{equation}
 	Using (\ref{Lemma_using_Tutte_Delta}), note that as $n \rightarrow \infty$ 
 	 \begin{equation}\label{Lemma_using_Tutte_DeltaB-L}
 	 		||B_{diag} - \Lambda||_2  = O(n^{-1/2+3\varepsilon}),
 	 \end{equation}
 	 where $\Lambda$ denotes the diagonal matrix whose diagonal elements are equal to the components of the vector $Q\vec{\theta}$.
 	 Since the algebraic connectivity 	$\lambda_1 \geq \sigma n$, we get that as $n \rightarrow \infty$
 	 \begin{equation}\label{Bdiag-lambda}
 	 		\left|tr\left((B_{diag} - \Lambda)\hat{Q}^{-1}\right)\right| \leq n ||B_{diag} - \Lambda||_2||\hat{Q}^{-1}||_2 = O(n^{-1/2+3\varepsilon}).
 	 \end{equation}
 	 Using (\ref{Bskewsmall}) and (\ref{Bdiag-lambda}), we obtain that as $n \rightarrow \infty$
 	 \begin{equation}\label{Lemma_using_Tutte_Phi}
 	 		tr(\Phi) = tr(B_{diag}\hat{Q}^{-1}) = tr(\Lambda\hat{Q}^{-1})+ O(n^{-1/2+3\varepsilon}) = \vec{\theta}^T Q \vec{\alpha} + O(n^{-1/2+3\varepsilon}),
 	 \end{equation} 
 	 where $\vec{\alpha}$ denotes the vector composed of the diagonal elements of the matrix $\hat{Q}^{-1}$.
 	 
 	 Using the property of the trace function $tr(XY) = tr(YX)$, we have that
 	 \begin{equation}\label{Lemma_using_Tutte_Phi2}
 	 		tr \left(\Phi^2\right) = tr\left((B_{skew}\hat{Q}^{-1})^2\right) + tr\left((B_{diag}\hat{Q}^{-1})^2\right) + 
						2\, tr\left(B_{skew}\hat{Q}^{-1}B_{diag}\hat{Q}^{-1}\right).
 	 \end{equation}
 	 Since  $B_{skew}$ is the skew-symmetric matrix and $\hat{Q}^{-1}B_{diag}\hat{Q}^{-1}$ is the symmetric matrix
 	 \begin{equation}\label{Lemma_using_Tutte_B0}
 	 		tr\left(B_{skew}\hat{Q}^{-1}B_{diag}\hat{Q}^{-1}\right) = 0.
 	 \end{equation}
 	 One can show that
 	 \begin{equation} \label{HS}
 	 	\begin{aligned}
 	 		tr \left( X^2\right) &\leq ||X||_{HS}^2, \\
 	 		||XY||_{HS} &\leq ||X||_{HS}||Y||_2.
 	 	\end{aligned}
 	 \end{equation}	
 	 	Therefore we get that
 	 \begin{equation}\label{Lemma_using_Tutte_trBsk}
 	 		\left|tr\left((B_{skew}\hat{Q}^{-1})^2\right)\right| \leq ||B_{skew}\hat{Q}^{-1}||^2_{HS}. 
 	 \end{equation}
 	 Since the algebraic connectivity 	$\lambda_1 \geq \sigma n$, using (\ref{Lemma_using_Tutte_Delta}), we obtain that 
 	 as $n \rightarrow \infty$ 
 	 \begin{equation}\label{Lemma_using_Tutte_Bsk}
 	 	||B_{skew}\hat{Q}^{-1}||_{HS} \leq ||\hat{Q}^{-1}||_2||B_{skew}||_{HS} 
 	 	\leq
 	 	 \frac{1}{\lambda_1}||B_{skew}||_{HS} = O(n^{-1/2 +\varepsilon}).
 	 \end{equation}
 	 Using (\ref{Lemma_using_Tutte_Delta}) and (\ref{Lemma_using_Tutte_DeltaB-L}), we get that
 	 \begin{equation}
 	 		\left|tr\left( (B_{diag}-\Lambda)\hat{Q}^{-1}B_{diag}\hat{Q}^{-1}\right)\right| \leq 
 	 		n \frac{1}{\lambda_1^2} ||(B_{diag}-\Lambda)||_2||B_{diag}||_2 = O(n^{-1+4\varepsilon})
 	 \end{equation}
 	 and
 	 \begin{equation}
 	 		\left|tr\left( (B_{diag}-\Lambda)\hat{Q}^{-1}(B_{diag}-\Lambda)\hat{Q}^{-1}\right)\right| \leq 
 	 		n \frac{1}{\lambda_1^2} ||(B_{diag}-\Lambda)||^2_2 = O(n^{-2+6\varepsilon}).
 	 \end{equation}
 	 Thus
 	 \begin{equation}\label{Lemma_using_Tutte_Bdiag2}
 	 		 tr\left((B_{diag}\hat{Q}^{-1})^2\right) =  tr\left((\Lambda\hat{Q}^{-1})^2\right) + O(n^{-1+4\varepsilon}).
 	 \end{equation}
 	 Combining (\ref{Lemma_using_Tutte_Phi2}), (\ref{Lemma_using_Tutte_B0}), (\ref{Lemma_using_Tutte_trBsk}), (\ref{Lemma_using_Tutte_Bsk}) 
 	 and (\ref{Lemma_using_Tutte_Bdiag2}), we obtain that
 	 \begin{equation}\label{Lemma_using_Tutte_Phi2f}
 	 		tr \left(\Phi^2\right) = tr\left((\Lambda\hat{Q}^{-1})^2\right) + O(n^{-1 +4\varepsilon}).
 	 \end{equation}
 	   Using (\ref{Lemma_using_Tutte_Phi}) and (\ref{Lemma_using_Tutte_Phi2f}) in (\ref{Lemma_using_Tutte_I+Phi}), we get that
 	  \begin{equation}\label{Lemma_using_Tutte_003}
 	  	\det(I+i \Phi) = \exp \left( i\,\vec{\theta}^T Q \vec{\alpha} + \frac{tr\left((\Lambda\hat{Q}^{-1})^2\right)}{2} + O(n^{-1/2+3\varepsilon}) \right).
 	  \end{equation}
 	  Combining  (\ref{Lemma_using_Tutte_000}), (\ref{Lemma_using_Tutte_001}) and (\ref{Lemma_using_Tutte_003}), we obtain (\ref{eq_Lemma_using_Tutte}).
\end{Proof}

We denote by $P(\vec{\theta})$ the orthogonal projection onto the space $L$. 
\begin{Lemma}\label{Lemma_Integral}
	Let the assumptions of Theorem \ref{main} hold. Let $\hat{Q} = Q + J$, where $J$ denotes the matrix with every entry 1. 
	For positive constants $a,b,c,d_1,d_2$ let sequence of vectors $\{\vec{\alpha}_n\}$ and 
	sequence of differentiable functions $\{R_n(\vec{\theta})\}$ be such that 
	\begin{equation}
		||\vec{\alpha}_n||_{\infty} \leq c/n,
	\end{equation}
	\begin{equation}
		|R_n(\vec{\theta})| \leq d_1\frac{\vec{\theta}^T \hat{Q}\vec{\theta}}{n},
	\end{equation}
	\begin{equation}
		R_n(\vec{\theta}) = R_n(P(\vec{\theta}))
	\end{equation}
	and for $\vec{\theta} \in U_n(\frac{4}{\sigma}n^{-1/2+\varepsilon})$
	\begin{equation}
	 \left\|\frac{\partial R_n(\vec{\theta})}{\partial\vec{\theta}}\right\|_{\infty} \leq d_2 \,n^{-1/2+\varepsilon}.
	\end{equation}
	For $n\geq 2$ define
	\begin{equation}\label{Integral_Lemma_Integral}
		J_n = \int\limits_{L \cap V_0}	
			\exp\left(
			i\,\vec{\theta}^T Q \vec{\alpha_n} - a\sum\limits_{(v_j,v_k)\in EG}\Delta_{jk}^2
			- b\sum\limits_{(v_j,v_k)\in EG} \Delta_{jk}^4 + R_n(\vec{\theta})
			\right)
		 dL.
	\end{equation}
	Then as $n\rightarrow\infty$ 
	\begin{equation}\label{eq_Lemma_Integral}
		J_n = \Theta_{k_1,k_2}
		\left(
			{\pi^{\frac{n-1}{2}} a^{-\frac{n-1}{2}} n^{1/2}  } \Big/ {\sqrt{\det\hat{Q}}}
		\right),
	\end{equation}
	where constants $k_1,k_2 > 0$ depend only on $a$, $b$, $c$, $d_1$, $d_2$ and $\sigma$.
\end{Lemma}
Lemma \ref{Lemma_Integral} is proved in Section 8.

\begin{Lemma}\label{Lemma_S0}
Let the assumptions of Theorem \ref{main} hold. Then as $n\rightarrow \infty$
\begin{equation}\label{eq_Lemma_S0}
	S_0 = 	\Theta_{k_1,k_2} 
		\left(
		2^{\frac{n-1}{2}}\pi^{\frac{n+1}{2}} n^{-1}\sqrt{\det\hat{Q}}
		\right),
\end{equation}
where constants $k_1,k_2 > 0$ depend only on $\sigma$.
\end{Lemma}
\begin{Proof} {\it Lemma \ref{Lemma_S0}.}
	Using formula (\ref{S_0}) with $r=1,2\ldots, n$ and summing, we obtain that
	\begin{equation}\label{Lemma_SO_nS0}
		\begin{aligned}
			nS_0 &= \sum\limits_{r=1}^{n}  \pi n^{1/2}\int\limits_{L \cap V_0}
    		 \prod\limits_{(v_j,v_k)\in EG} \cos \Delta_{jk} \sum \limits_{T \in {\cal T}_r} 
     		\prod \limits_{(v_j,v_k)\in ET} (1+i\tan \Delta_{jk})   
  	\ d L=\\
  		&= 
  			\pi n^{1/2}
  		  \int\limits_{L \cap V_0}
    		 \prod\limits_{(v_j,v_k)\in EG} \cos \Delta_{jk} \sum\limits_{r=1}^{n}\sum \limits_{T \in {\cal T}_r} 
     		\prod \limits_{(v_j,v_k)\in ET} (1+i\tan \Delta_{jk})   
  	\ 	dL.
  	\end{aligned}
	\end{equation}
	By Taylor's theorem  we have that for $\vec{\theta} \in V_0$
	\begin{equation}\label{Lemma_SO_Taylor}
		\prod\limits_{(v_j,v_k)\in EG} \cos \Delta_{jk} = 
		\exp\left(-\frac{1}{2}\sum\limits_{(v_j,v_k)\in EG}\Delta_{jk}^2 - 
			\frac{1}{12}\sum\limits_{(v_j,v_k)\in EG}\Delta_{jk}^4 + O(n^{-1+6\varepsilon})\right).
	\end{equation}
	Combining (\ref{Lemma_SO_nS0}) with (\ref{Lemma_SO_Taylor}) and Lemma \ref{Lemma_using_Tutte}, we obtain that as $n\rightarrow \infty$
	\begin{equation}\label{Lemma_S0_S0}
		S_0 \sim \pi\frac{\det\hat{Q}}{n^{3/2}}
		\int\limits_{L \cap V_0}	
		\exp\left(
				i\,\vec{\theta}^T Q \vec{\alpha} - \frac{1}{2}\sum\limits_{(v_j,v_k)\in EG}\Delta_{jk}^2
			- \frac{1}{12}\sum\limits_{(v_j,v_k)\in EG} \Delta_{jk}^4 + \frac{1}{2}R(\vec{\theta})
			\right) 
		 dL,
	\end{equation}
	where $\vec{\alpha}$ denotes the vector composed of the diagonal elements of the matrix $\hat{Q}^{-1}$ and
	\begin{equation}
		R(\vec{\theta}) = tr(\Lambda\hat{Q}^{-1})^2,
	\end{equation}
	where $\Lambda$ is the diagonal matrix whose diagonal elements are equal to the components of the vector $Q\vec{\theta}$. Since vector $[1, 1 \ldots, 1]^T$
	is the eigenvector of $Q$, corresponding to eigenvalue $0$, we have that	
	\begin{equation}
		Q\vec{\theta} = QP(\vec{\theta}).
	\end{equation}
	Thus
	\begin{equation}\label{PR}
		R(\vec{\theta}) = R(P(\vec{\theta})). 
	\end{equation}
	Note that Lemma \ref{Lemma_minor} implies as $n \rightarrow \infty$ 
	\begin{equation}\label{maxalpha}
		||\vec{\alpha}||_\infty \leq c_1/n,
	\end{equation}
	where $c_1 = c_1(\sigma) > 0$. Since the algebraic connectivity 	$\lambda_1 \geq \sigma n$,	using (\ref{HS}), we get that
	\begin{equation}\label{maxR}
		\begin{aligned}
				|R(\vec{\theta})| &\leq ||\Lambda\hat{Q}^{-1}||_{HS}^2 \leq ||\Lambda||_{HS}^2 ||\hat{Q}^{-1}||_2^2 = \\
		 		&=||Q\vec{\theta}||_2^2\ ||\hat{Q}^{-1}||_2^2 \leq ||Q||_2^2\ ||\hat{Q}^{-1}||_2^2\ ||\vec{\theta}||_2^2 \leq \\
		 		&\leq \lambda_{n-1}^2 \frac{1}{\lambda_1^2} ||\vec{\theta}||_2^2 \leq   n^2 \frac{1}{\lambda_1^2} \frac{\vec{\theta}^T \hat{Q}\vec{\theta}}{\lambda_1}
		 		\leq \frac{1}{\sigma^3}\frac{\vec{\theta}^T \hat{Q}\vec{\theta}}{n}.
		\end{aligned}
	\end{equation}
	Note that for $\vec{\theta} \in U_n(\frac{4}{\sigma}n^{-1/2+\varepsilon})$ and for some $1\leq k \leq n$
	\begin{equation}\label{Lambda2}
		||\Lambda||_2 \leq \sum\limits_{j=1}\limits^{n} |\Delta_{jk}| = O\left( n^{1/2+\varepsilon} \right).
	\end{equation}
	For $1 \leq k \leq n$ we denote by $(\hat{Q}^{-1})_k$ the $k$-th column of the matrix $\hat{Q}^{-1}$.
	Using again $\lambda_1 \geq \sigma n$, we get that 
	\begin{equation}\label{Q-1k}
		1 = ||\hat{Q} (\hat{Q}^{-1})_k ||_2 \geq \lambda_1 ||(\hat{Q}^{-1})_k||_2  \geq \sigma n ||(\hat{Q}^{-1})_k||_2. 
	\end{equation}
	Note that
	\begin{equation}
		\frac{\partial R}{\partial \theta_k} = 
 		2tr(\frac{\partial \Lambda}{\partial \theta_k}\hat{Q}^{-1}\Lambda\hat{Q}^{-1})
 		= 2d_k (\hat{Q}^{-1}\Lambda\hat{Q}^{-1})_{kk} + 2tr(\tilde{\Lambda}\hat{Q}^{-1}\Lambda\hat{Q}^{-1}), 
	\end{equation}
	where
		$(\hat{Q}^{-1}\Lambda\hat{Q}^{-1})_{kk}$ denotes  $(k,k)$-th element of the matrix $\hat{Q}^{-1}\Lambda\hat{Q}^{-1}$ 
	and the matrix $\tilde{\Lambda}$ is such that for any $1 \leq j \leq n$ the diagonal element $|\Lambda_{jj}|\leq 1$  
.
	Since the algebraic connectivity 	$\lambda_1 \geq \sigma n$, using (\ref{Lambda2}) (\ref{Q-1k}), we get that as $n\rightarrow \infty$
	\begin{equation}
		|tr(\tilde{\Lambda}\hat{Q}^{-1}\Lambda\hat{Q}^{-1})| \leq n \, ||\tilde{\Lambda}\hat{Q}^{-1}\Lambda\hat{Q}^{-1}||_2 
		\leq n \, ||\tilde{\Lambda}||_2 \  ||\Lambda||_2 \ \frac{1}{\lambda_1^2} 
		= O(n^{-1/2+\varepsilon})
	\end{equation}
	and 
	\begin{equation}
		2d_k (\hat{Q}^{-1}\Lambda\hat{Q}^{-1})_{kk} = 2d_k \Lambda_{kk} ||(\hat{Q}^{-1})_k||_2^2 = O(n^{-1/2+\varepsilon}).
	\end{equation}
	Thus for some $d>0$, depending only on $\sigma$
	\begin{equation}\label{partialR}
		\left\|\frac{\partial R_n(\vec{\theta})}{\partial\vec{\theta}}\right\|_{\infty} \leq d\, n^{-1/2+\varepsilon}.
	\end{equation}
	Combining (\ref{Lemma_S0_S0}), (\ref{PR}), (\ref{maxalpha}), (\ref{maxR}), (\ref{partialR}) and using Lemma \ref{Lemma_Integral} we obtain (\ref{eq_Lemma_S0}).
\end{Proof}
 
%\begin{equation}\label{S_0} 
%S_0 = \pi n^{1/2} \int_{L \cap U_{n}(n^{-1/2 + \varepsilon})}
%\left( 
%     \prod\limits_{(v_j,v_k)\in EG} \cos \Delta_{jk} \sum \limits_{T \in {\cal T}_r} 
%     \prod \limits_{(v_j,v_k)\in ET} (1+i\tan \Delta_{jk})  
%  \right) 
%  dL,
%\end{equation}
%%%%%%%%%%%%%%%%%%%%%%%%%%%%%%%%%%%%%%%%%%%%%%%%%%%%%%%%%%%%%%%%%%%%%%%%%%%%%%%%%%%%%%%%%%%%%%%%%%%%%%%%%%%%%%%%%%%%%%%%%
%%%%%%%%%%%%%%%%%%%%%%%%%%%%%%%%%%%%%%%%%%%%%%%%%%%%%%%%%%%%%%%%%%%%%%%%%%%%%%%%%%%%%%%%%%%%%%%%%%%%%%%%%%%%%%%%%%%%%%%%%
%%%%%%%%%%%%%%%%%%%%%%%%%%%%%%%%%%%%%%%%%%%%%%%%%%%%%%%%%%%%%%%%%%%%%%%%%%%%%%%%%%%%%%%%%%%%%%%%%%%%%%%%%%%%%%%%%%%%%%%%%
%%%%%%%%%%%%%%%%%%%%%%%%%%%%%%%%%%%%%%%%%%%%%%%%%%%%%%%%%%%%%%%%%%%%%%%%%%%%%%%%%%%%%%%%%%%%%%%%%%%%%%%%%%%%%%%%%%%%%%%%%

\section{The insignificant parts of the integral}
In this section we  prove that $S_0$ contributes almost all of $S$, even though it involves
only a tiny part of the region of integration, compare with Section 4 of \cite{Brendan1995}. 
We continue to use the same value of $\varepsilon$ as in
the previous section. 

Let assumptions of Theorem 2.1 hold.
Define $E'T = \{(v_j,v_k), (v_k,v_j) \ | \ (v_j,v_k)\in ET \}$.  
We express the integrand of (\ref{S}) as
\begin{equation}
	F(\vec{\theta}) = \sum\limits_{T\in T_r} \prod\limits_{jk \in EG} f_{jk}(T,\vec{\theta}),
\end{equation}
where
\begin{equation}
		f_{jk}(T,\vec{\theta}) = 
		\left\{
			\begin{array}{ll}
				\cos \Delta_{jk} (1+i\tan \Delta_{jk}),& (v_j,v_k)\in ET,\\
				\cos \Delta_{jk} (1-i\tan \Delta_{jk}),& (v_k,v_j)\in ET,\\
				\cos \Delta_{jk},& \text{otherwise.}
			\end{array}
\right.
\end{equation}
Note that $|f_{jk}(T,\vec{\theta})|\leq 1$ for all values of the parameters. One can show that
\begin{equation}\label{cosx}
	|\cos(x)|\leq \exp (-\frac{1}{2} x^2) \ \text{ for }  \ |x|\leq \frac{9}{16}\pi.
\end{equation}
	
	Divide the interval $[-\frac{1}{2}\pi,\frac{1}{2}\pi] \text{ mod } \pi$ into 32 equal intervals $H_0,\ldots,H_{31}$ such that $H_0 = [-\frac{1}{64}\pi,\frac{1}{64}\pi]$. For each $j$, define the region $W_j\subseteq U_n(\pi/2)$ as the set of points having at least $\frac{1}{32}n$ 
	coordinates in $H_j$. Clearly, the $W_j$'s cover $U_n(\pi/2)$ and also each $W_j$
can be mapped to $W_0$ by a uniform translation of the $\theta_j \text{ mod } \pi$. This mapping preserves the integrand of (\ref{S}) and also maps 
$V_0$ to itself, so we have that
\begin{equation}
	\int\limits_{U_n(\pi/2)-V_0}|F(\theta)|d\vec{\theta}\leq 32Z,
\end{equation}
where
\begin{equation}
	Z=\int\limits_{W_0-V_0}|F(\theta)|d\vec{\theta}.
\end{equation}

	We proceed by defining integrals $S_1,\ldots,S_4$ in such a way that $Z$ is obviously bounded by their sum. We then show 
that $S_j = o(S_0)$ for $j=1,2,3,4$ separately. Write 
\begin{equation}
	F(\vec{\theta}) = F_a(\vec{\theta}) + F_b(\vec{\theta}),
\end{equation}
where $F_a(\theta)$ and $F_b(\theta)$ are defined by restricting the sum to trees with maximum
degree greater than $\sigma n/4$ and no more than $\sigma n/4$, respectively. Also define regions $V_1$ and $V_2$ as
follows.
\begin{equation}
	\begin{aligned}
	&V_1 =\{\vec{\theta}\in W_0 \text{ $|$ } |\theta_j|\geq \frac{1}{32}\pi \text{ for fewer than $n^\varepsilon$ values of $j$} \},\\
	&V_2 =\{\vec{\theta}\in V_1 \text{ $|$ } |\theta_j|\geq \frac{1}{16}\pi \text{ for at least one value of $j$} \}.
	\end{aligned}
\end{equation}
Then our four integrals can be defined as
\begin{equation}
	\begin{aligned}
	&S_1 = \int\limits_{W_0-V_1}|F(\vec{\theta})|d\vec{\theta},\\
	&S_2 = \int\limits_{V_1}|F_a(\vec{\theta})|d\vec{\theta},\\
	&S_3 = \int\limits_{V_2}|F_b(\vec{\theta})|d\vec{\theta},\\
	&S_4 = \int\limits_{V_1-V_2-V_0}|F_b(\vec{\theta})|d\vec{\theta}.
	\end{aligned}
\end{equation}

We begin with $S_1$. Let $h$ be the function from Lemma \ref{Lemma_struct} for the set	$A = \{v_j \ | \ |\theta_j|\leq\frac{1}{64}\pi\}$.
	We denote $l_{min}$ such natural number that inequality
	\begin{equation}
		|\theta_j| \geq \frac{1}{64}\pi (1 + l/H)
	\end{equation} 
	holds for at least $n^\varepsilon/H$ indices of the set $\{j \ | \ h(v_j) = l\}$. Existence of $l_{min}$ follows from the definition of the
	region $V_1$. 
	If $\theta_j$ and $\theta_k$ are such that
	\begin{equation} 
		|\theta_j|\geq \frac{1}{64}\pi (1 + l_{min}/H) \  \text{ and } \  |\theta_k|\leq \frac{1}{64}\pi (1 + (l_{min}-1)/H)
	\end{equation} 
	or vice versa, but 
$(v_j,v_k)\notin E'T$, we have that
	$|f_{jk}(T,\vec{\theta})|\leq \cos(\frac{1}{64}\pi/H).$ This includes at least  
$ (\alpha n - n^{\varepsilon} ) {\frac{n^\varepsilon}{H}} - n$ edges $(v_j,v_k)\in EG$. Using (\ref{lambda<n}) and (\ref{t(G)hat}), we get that as $n\rightarrow \infty$
\begin{equation}\label{S_1}
	S_1\leq t(G)\pi^n\left(\cos\frac{\pi}{64H}\right)^{(\alpha n - n^{\varepsilon} ) \frac{n^\varepsilon}{H}-n} = 
	O\Big(\exp(-cn^{1+\varepsilon}) \Big)2^{\frac{n-1}{2}}\pi^{\frac{n+1}{2}} n^{-1}\sqrt{\det\hat{Q}}
\end{equation}
for some constant $c > 0$ depending only on $\sigma$.

To bound $S_2$, we first note from Lemma \ref{Lemma_Trees>d} that the number of trees with maximum de-
gree greater than $\sigma n /4$ is less than $c_3^n \det\hat{Q}/(\sigma n /4)!$. Using (\ref{cosx}), we see that
\begin{equation}\label{f_jkexp}
	|f_{jk}(T,\vec{\theta})|\leq \exp(-\frac{1}{2}\Delta_{jk}^2)
\end{equation}
except for at most $n^{2\varepsilon}$ pairs $(j,k)$
with $|\Delta_{jk}|\geq\frac{1}{16}\pi$ and fewer than $n$ pairs in $E'T$.
In those excluded cases the value $\exp(-\frac{1}{2}\Delta_{jk}^2)$ may be high by a factor $\exp(\frac{1}{2}\pi^2)$.
Hence, we have that
\begin{equation}\label{S_21}
 S_2\leq \frac{c_3^n \det\hat{Q}}{(\sigma n /4)!} 
 \exp\left(\frac{1}{2}\pi^2(n+n^{2\varepsilon})\right)
 \int\limits_{U_n(\pi/2)}\exp\left(-\frac{1}{2}\sum\limits_{(v_j,v_k)\in EG}\Delta_{jk}^2\right)d\vec{\theta}. 
\end{equation}
\begin{Lemma}\label{Lemma_expUn}
	Let the assumptions of Lemma \ref{Lemma_normQ} hold. Then 
	\begin{equation}\label{eq_Lemma_expUn}
		\int\limits_{U_n(\pi/2)}\exp\left(-\frac{1}{2}\sum\limits_{(v_j,v_k)\in EG}\Delta_{jk}^2\right)d\vec{\theta}
		\leq \frac{ 2^{\frac{n-1}{2}} \pi^{\frac{n+1}{2}} n}{\sqrt{\det{\hat{Q}}}}.
	\end{equation}
\end{Lemma}

\noindent
Lemma \ref{Lemma_expUn} is proved in Section 8. Combining (\ref{S_21}) and (\ref{eq_Lemma_expUn}),  we obtain that as $n\rightarrow \infty$
\begin{equation}\label{S_2}
	S_2 = O(n^{-cn})\, 2^{\frac{n-1}{2}}\pi^{\frac{n+1}{2}} n^{-1}\sqrt{\det\hat{Q}}
\end{equation}
for some constant $c > 0$ depending only on $\sigma$.

We denote by $G_T$ the graph which arises
from $G$ by removing all edges of the tree $T$. Let $G_{T,\vec{\theta}}$ be the graph resulting from
from $G_T$ by removing vertices, corresponding to those values of $j$ for which $|\theta_j| \geq \frac{1}{16}\pi$.

For $1\leq r\leq n^\varepsilon$ let $S_3(r)$ denote the contribution to $S_3$ 
of those $\theta \in V_2$ such that 
$|\theta_j|\geq\frac{1}{16}\pi$ for exactly $r$ values of $j$. 
If $|\theta_j|\leq\frac{1}{32}\pi$ and $|\theta_k|\geq\frac{1}{16}\pi$ 
or vice versa, we have that 
\begin{equation}\label{f_jkcos}
	|f_{jk}(T,\vec{\theta})|\leq\cos\left(\frac{1}{32}\pi\right)
\end{equation}
unless $(v_j,v_k) \in E'T$.
This includes at least $r(\sigma n/2 -\sigma n/4 - n^\varepsilon)$ pairs $(j,k)$, because 
 the degree of any vertex of the graph $G$ is at least $\sigma n /2$, see (\ref{lambda<d}).
For pairs $(j,k)$ such that $|\theta_j|,|\theta_k|\leq\frac{1}{16}\pi$, but $(v_j,v_k)\notin E'T$, we use (\ref{f_jkexp}).
We put $\vec{\theta'} = (\theta_1,\ldots,\theta_{n-r})$.
Then,
\begin{equation}\label{S_3GTtheta}
	S_3(r)\leq \pi^r  \left(\cos\frac{\pi}{32}\right)^{r(\sigma n/4 - n^\varepsilon)} 
	\sum\limits_{(r)} \sum\limits_T \int\limits_{ U_{n-r}(\pi/2) }
	\exp\left(-\frac{1}{2}\sum\limits_{(v_j,v_k)\in EG_{T,\vec{\theta}}}\Delta_{jk}^2\right)
	d\vec{\theta'},
\end{equation}
 where the first sum is over choices of those values of $j$ for which $|\theta_j| \geq \frac{1}{16}\pi$ and
 the second sum is over trees with maximum degree $\sigma n /4$. 
 Using Lemma \ref{Lemma_G_T} and then Lemma \ref{Lemma_connectivity} and Lemma \ref{Lemma_Gr} for the graph $G_T$, we obtain that 
 \begin{equation}
 	\lambda_1(G_{T,\vec{\theta}}) \geq \sigma n /2 -  n^\varepsilon 
 \end{equation}	
 	and 
 \begin{equation}\label{GTtheta}
 		\det{\hat{Q}(G_{T,\vec{\theta}})} \geq \frac{\det\hat{Q}}{\left(c_5 n\right)^r},
 \end{equation}
 where $c_5 = c_5(\sigma) > 0$ and $\hat{Q} = \hat{Q}(G)$. According to Lemma \ref{Lemma_expUn}, we have that
 \begin{equation}\label{IntegralGTtheta}
 		\int\limits_{U_{n-r}(\pi/2)}\exp\left(-\frac{1}{2}\sum\limits_{(v_j,v_k)\in EG_{T,\vec{\theta}}}\Delta_{jk}^2\right)d\vec{\theta'}
		\leq \frac{ 2^{\frac{n-r-1}{2}} \pi^{\frac{n-r+1}{2}} n}{\sqrt{\det{\hat{Q}(G_{T,\vec{\theta}})}}}.
 \end{equation}
 Combining (\ref{S_3GTtheta}) with  (\ref{GTtheta}), (\ref{IntegralGTtheta}) and allowing $n^r$ for the choice 
 of those values of $j$ for which $|\theta_j| \geq \frac{1}{16}\pi$, we obtain that
 \begin{equation}
  	 		S_3(r) \leq  2^{\frac{n-r-1}{2}} \pi^{\frac{n+r+1}{2}} n^{r+1} 
 	 		\left(\cos\frac{\pi}{32}\right)^{r(\sigma n/4 - n^\varepsilon)}  \frac{t(G)(c_5 n)^{r/2}}{\sqrt{\det\hat{Q}}} 
  \end{equation} 
  and, using (\ref{t(G)hat}), we can calculate that 
 \begin{equation}\label{S_3}
 		S_3 = \sum\limits_{r=1}\limits^{n^\varepsilon}S_3(r) = O(c^{-n})\, 2^{\frac{n-1}{2}}\pi^{\frac{n+1}{2}} n^{-1}\sqrt{\det\hat{Q}}
 \end{equation}
 for some constant $c>1$ depending only on $\sigma$.
 
 		Since $\Delta_{jk}\leq\frac{1}{8}\pi$ for $\vec{\theta} \in V_1 - V_2 - V_0$ and the integrand is invariant under uniform translation of all
	the $\theta_j$'s mod $\pi$, we can fix $\bar{\theta} = 0$ and multiply it by the ratio 
	of its range $\pi$ to the length $n^{-1/2}$ of the vector $\frac{1}{n}[1,1,\ldots,1]^T$. 
	Thus we get that
	\begin{equation} 
		S_4 \leq \pi n^{1/2} \int\limits_{L\cap U_n({\pi}/{8}) - V_0}
	|F_b(\vec{\theta})|dL,
	\end{equation}
where $L$ denotes the orthogonal complement to the vector $[1,1,\ldots,1]^T$.
 In a similar way as (\ref{S_3GTtheta}) we find that
 	\begin{equation}\label{S44}
 		S_4\leq  	\pi n^{1/2}	\sum\limits_T \int\limits_{ L\cap U_n({\pi}/{8}) - V_0}
	\exp\left(-\frac{1}{2}\sum\limits_{(v_j,v_k)\in EG_T}\Delta_{jk}^2\right)
	dL,
 	\end{equation}
 where the first sum is over trees with maximum degree $\sigma n /4$. 
 \begin{Lemma}\label{Lemma_exp-V0}
	Let the assumptions of Lemma \ref{Lemma_normQ} hold. Then as $n\rightarrow \infty$
	\begin{equation}\label{eq_Lemma_exp-V0}
		\int\limits_{L -U_n(n^{-1/2+\varepsilon})}\exp\left(-\frac{1}{2}\sum\limits_{(v_j,v_k)\in EG}\Delta_{jk}^2\right)dL
		=O\left(\exp(-cn^{2\varepsilon})\right) \frac{ 2^{\frac{n-1}{2}} \pi^{\frac{n-1}{2}} n^{1/2}}{\sqrt{\det{\hat{Q}}}}
	\end{equation}
	for some $c>0$ depending only on $\sigma$.
\end{Lemma}

\noindent
Lemma \ref{Lemma_exp-V0} is proved in Section 8. Using Lemma \ref{Lemma_G_T} and  combining (\ref{eq_Lemma_exp-V0}), (\ref{S44}) and (\ref{t(G)hat}),
we obtain that as $n\rightarrow \infty$
\begin{equation}\label{S_4}
	S_4 = O\left(\exp(-cn^{2\varepsilon})\right) t( G)\frac{ 2^{\frac{n-1}{2}} \pi^{\frac{n+1}{2}} n}{\sqrt{\det{\hat{Q}}}}
		= O\left(\exp(-cn^{2\varepsilon})\right) 2^{\frac{n-1}{2}} \pi^{\frac{n+1}{2}} n^{-1}\sqrt{\det{\hat{Q}}}
\end{equation}
	for some $c>0$ depending only on $\sigma$.
	Combining (\ref{S_1}), (\ref{S_2}), (\ref{S_3}), (\ref{S_4}) and Lemma \ref{Lemma_S0}, we obtain the desired result.
\begin{Lemma}\label{Lemma_S1234}
	Let the assumptions of Theorem \ref{main} hold. Then as $n\rightarrow \infty$
	\begin{equation}
		S  = \left(1 + O\left(\exp(-cn^{2\varepsilon})\right)\right) S_0
	\end{equation}
	for some $c>0$ depending only on $\sigma$.
\end{Lemma}
%%%%%%%%%%%%%%%%%%%%%%%%%%%%%%%%%%%%%%%%%%%%%%%%%%%%%%%%%%%%%%%%%%%%%%%%%%%%%%%%%%%%%%%%%%%%%%%%%%%%%%%%%%%%%%%%%%%%%%%%%
%%%%%%%%%%%%%%%%%%%%%%%%%%%%%%%%%%%%%%%%%%%%%%%%%%%%%%%%%%%%%%%%%%%%%%%%%%%%%%%%%%%%%%%%%%%%%%%%%%%%%%%%%%%%%%%%%%%%%%%%%
%%%%%%%%%%%%%%%%%%%%%%%%%%%%%%%%%%%%%%%%%%%%%%%%%%%%%%%%%%%%%%%%%%%%%%%%%%%%%%%%%%%%%%%%%%%%%%%%%%%%%%%%%%%%%%%%%%%%%%%%%
%%%%%%%%%%%%%%%%%%%%%%%%%%%%%%%%%%%%%%%%%%%%%%%%%%%%%%%%%%%%%%%%%%%%%%%%%%%%%%%%%%%%%%%%%%%%%%%%%%%%%%%%%%%%%%%%%%%%%%%%%

\section{Proof of Lemma \ref{Lemma_normQ}}
According to (\ref{EulS}) and (\ref{S})

\begin{equation}
Eul(G) = \prod\limits_{j=1}^n \big( \frac{d_j}{2} -1 \big)! \, 2^{|EG|-n+1} \pi^{-n} S,
\end{equation}
where
\begin{equation}
	S = \int\limits_{U_n(\pi/2)} 
 % \left( 
     \prod\limits_{(v_j,v_k)\in EG} \cos \Delta_{jk} \sum \limits_{T \in {\cal T}_r} 
     \prod \limits_{(v_j,v_k)\in ET} (1+i\tan \Delta_{jk})  
 % \right) 
 \
 d\vec{\theta}.
\end{equation}
Combining Lemma \ref{Lemma_S0} and Lemma \ref{Lemma_S1234} we get that as $n\rightarrow \infty$
\begin{equation}
		S = 	\Theta_{k_1,k_2} 
		\left(
		2^{\frac{n-1}{2}}\pi^{\frac{n+1}{2}} n^{-1}\sqrt{\det\hat{Q}}
		\right),
\end{equation}
where constants $k_1,k_2 > 0$ depend only on $\sigma$. Taking into account (\ref{t(G)hat}) we obtain (\ref{main_eq}).
{\scriptsize $\blacksquare$}

\vspace{3mm}

If for some $\sigma> 1/2$ the degree of each vertex of the graph $G$ at least $\sigma n$, we can use (\ref{lambda>d}) and get that	
\begin{equation}
	\lambda_1(G) \geq 2 \min_j d_j -n + 2 > (2\sigma - 1)n.
\end{equation}

%\noindent
%{\bf Remark 7.1.} 
%	In fact, using Lemma \ref{Lemma_using_Tutte} and Lemma \ref{Lemma_S1234}, the estimation of the number of the Eulerian circuits
%	 is reduced (see proof of Lemma \ref{Lemma_S0}) to estimating the integral
%\begin{equation}\label{888}
%		\int\limits_{V_0}	
%		\exp\left(
%				i\,\vec{\theta}^T Q \vec{\alpha} - \frac{1}{2}\sum\limits_{(v_j,v_k)\in EG}\Delta_{ij}^2
%			- \frac{1}{12}\sum\limits_{(v_j,v_k)\in EG} \Delta_{jk}^4 + tr(\Lambda \hat{Q}^{-1})^2
%			\right) 
%		 d\vec{\theta},
%\end{equation}
%where $\alpha$ denotes the vector composed of the diagonal elements of $\hat{Q}^{-1}$, $\Lambda$ denotes
%the diagonal matrix whose diagonal elements  are equal to components of the vector $Q\vec{\theta}$. Apparently, it is possible 
%to estimate integral (\ref{888})
% more accurately for particular classes of graphs and obtain  
%asymptotic formulas for $Eul(G)$, similar to (\ref{Brendanf}).

%%%%%%%%%%%%%%%%%%%%%%%%%%%%%%%%%%%%%%%%%%%%%%%%%%%%%%%%%%%%%%%%%%%%%%%%%%%%%%%%%%%%%%%%%%%%%%%%%%%%%%%%%%%%%%%%%%%%%%%%%
%%%%%%%%%%%%%%%%%%%%%%%%%%%%%%%%%%%%%%%%%%%%%%%%%%%%%%%%%%%%%%%%%%%%%%%%%%%%%%%%%%%%%%%%%%%%%%%%%%%%%%%%%%%%%%%%%%%%%%%%%
%%%%%%%%%%%%%%%%%%%%%%%%%%%%%%%%%%%%%%%%%%%%%%%%%%%%%%%%%%%%%%%%%%%%%%%%%%%%%%%%%%%%%%%%%%%%%%%%%%%%%%%%%%%%%%%%%%%%%%%%%
%%%%%%%%%%%%%%%%%%%%%%%%%%%%%%%%%%%%%%%%%%%%%%%%%%%%%%%%%%%%%%%%%%%%%%%%%%%%%%%%%%%%%%%%%%%%%%%%%%%%%%%%%%%%%%%%%%%%%%%%%

\section{Proofs of Lemma \ref{Lemma_Integral}, Lemma \ref{Lemma_expUn} and Lemma \ref{Lemma_exp-V0}}

Let assumptions of Lemma \ref{Lemma_normQ} hold. We define 
\begin{equation}
	\vec{\phi} = \vec{\phi}(\vec{\theta}) = (\phi_1(\vec{\theta}),\ldots \phi_n(\vec{\theta}))
	= \hat{Q}\vec{\theta}. 
\end{equation}
We continue to use notation $P(\vec{\theta)}$ for the orthogonal projection onto the space $L$, where 
$L$ is the orthogonal complement to the vector $[1,1,\ldots,1]^T$. For any $a>0$ we have that
\begin{equation}\label{Integral_hatQ}
	\int\limits_{\mathbb{R}^n} e^{-a\,\vec{\theta}^T\hat{Q}\vec{\theta}} d \vec{\theta} =  \pi^{{n}/{2}} a^{-{n}/{2}} \Big/ \sqrt{\det\hat{Q}}
\end{equation}
and
\begin{equation}\label{Integral_hatQL}
	\int\limits_{L} e^{-a\,\vec{\theta}^T\hat{Q}\vec{\theta}} d L
	=
	\int\limits_{L} e^{-a\,\vec{\theta}^TQ\vec{\theta}} d L =  \pi^{\frac{n-1}{2}} a^{-\frac{n-1}{2}} n^{1/2} \Big/ \sqrt{\det\hat{Q}}.
\end{equation}

\begin{Proof} {\it Lemma \ref{Lemma_expUn}.}
	Note that
	\begin{equation}\label{sum->Q}
		\sum\limits_{(v_j,v_k)\in EG}\Delta_{jk}^2= \vec{\theta}^T Q \vec{\theta}.
	\end{equation}
	Since diagonal of $U_n(\pi/2)$ is equal to $\pi n^{1/2}$ and $Q\vec{\theta} = QP(\vec{\theta})$ we have that
	\begin{equation}\label{Integral<Un}
		\int\limits_{U_n(\pi/2)}\exp\left(-\frac{1}{2}\sum\limits_{(v_j,v_k)\in EG}\Delta_{jk}^2\right)d\vec{\theta}
		\leq
		\pi n^{1/2} \int\limits_{L} e^{-\frac{1}{2}\,\vec{\theta}^T Q\vec{\theta}} d L.
	\end{equation}
	Using (\ref{Integral_hatQL}), we obtain (\ref{eq_Lemma_expUn}).
\end{Proof}

Note that for some $g_1(\vec{\theta}) = g_1(\theta_2,\ldots,\theta_{n})$ 
\begin{equation}\label{g_theta}
	\vec{\theta}^T\hat{Q}\vec{\theta} = \frac{\phi_1(\vec{\theta})^2}{d_1+1} + g_1(\vec{\theta}).
\end{equation}
Using (\ref{Lemma_normQ_d1}), we get that as $n\rightarrow \infty$
\begin{equation}
	\begin{aligned}
	\int\limits_{\mathbb{R}^n} e^{-a\,\vec{\theta}^T\hat{Q}\vec{\theta}} d \vec{\theta}
	&=
	\int\limits_{-\infty}\limits^{+\infty}\cdots \int\limits_{-\infty}\limits^{+\infty}
	e^{-a\,g_1(\theta_2,\ldots,\theta_{n})}
	\left( \int \limits_{-\infty}\limits^{+\infty} e^{-a\,\frac{\phi_1(\vec{\theta})^2}{d_1+1}} d \theta_1  \right)
	d \theta_2 \ldots d \theta_n\\
	&=
		\left(1 + O\left(\exp(-\tilde{c}n^{2\varepsilon})\right)\right)
			\int \limits_{|\phi_1(\vec{\theta})|\leq \frac{1}{2}c_\infty^{-1} n^{1/2+\varepsilon}} e^{-a\,\vec{\theta}^T\hat{Q}\vec{\theta}} d \vec{\theta}  
	\end{aligned}
\end{equation}
for some $\tilde{c}>0$ depending only on $\sigma$ and $a$, where $c_\infty$ is the constant of Lemma \ref{Lemma_normQ}.
Combining similar expressions for $\phi_1, \phi_2, \ldots \phi_n$, we obtain that as $n\rightarrow \infty$
\begin{equation}
	\int \limits_{||\vec{\phi}(\vec{\theta})||_\infty\leq \frac{1}{2}c_\infty^{-1} n^{1/2+\varepsilon}} e^{-a\,\vec{\theta}^T\hat{Q}\vec{\theta}} d \vec{\theta}
	= \left(1 + O\left(\exp(-cn^{2\varepsilon})\right)\right) \int\limits_{\mathbb{R}^n} e^{-a\,\vec{\theta}^T\hat{Q}\vec{\theta}} d \vec{\theta}
\end{equation}
for some $c>0$ depending only on $\sigma$ and $a$. Using Lemma \ref{Lemma_normQ}, we get that as $n\rightarrow \infty$
\begin{equation}\label{1/2U_n}
		\int \limits_{U_n( \frac{1}{2}n^{-1/2+\varepsilon})} e^{-a\,\vec{\theta}^T\hat{Q}\vec{\theta}} d \vec{\theta}
	= \left(1 + O\left(\exp(-cn^{2\varepsilon})\right)\right) \int\limits_{\mathbb{R}^n} e^{-a\,\vec{\theta}^T\hat{Q}\vec{\theta}} d \vec{\theta}.
\end{equation}

\begin{Proof} {\it Lemma \ref{Lemma_exp-V0}.}
	Note that 
	\begin{equation}
		||P(\vec{\theta})||_\infty = ||\vec{\theta} 
		- \bar{\theta} [1,1,\ldots,1]^T||_\infty \leq 2 ||\vec{\theta} ||_\infty,
	\end{equation}
	where
	\begin{equation}
 		\bar{\theta} = \frac{\theta_1 + \theta_2 + \ldots  \theta_n}{n}.
	\end{equation}
	Thus
	\begin{equation}\label{U_n_in_P}
		U_n( \frac{1}{2}n^{-1/2+\varepsilon}) \subset \left\{ \vec{\theta} \ | \ P(\vec{\theta}) \in U_n(n^{-1/2+\varepsilon})\right\}
	\end{equation}
	Since $Q\vec{\theta} = QP(\vec{\theta})$, using (\ref{sum->Q}) and (\ref{U_n_in_P}), we get that
	\begin{equation}\label{LcapU_n->Un}
		\begin{aligned}
			\int\limits_{L \cap U_n(n^{-1/2+\varepsilon})}\exp\left(-\frac{1}{2}\sum\limits_{(v_j,v_k)\in EG}\Delta_{jk}^2\right)dL =
		 	\int\limits_{L \cap U_n(n^{-1/2+\varepsilon})} e^{-\frac{1}{2}\,\vec{\theta}^T Q\vec{\theta}} d L \\
		 	= \int\limits_{P(\vec{\theta}) \in U_n(n^{-1/2+\varepsilon})} e^{-\frac{1}{2}\,\vec{\theta}^T \hat{Q}\vec{\theta}} d \vec{\theta}
		 	\Big/ \int\limits_{-\infty}\limits^{+\infty} e^{-\frac{1}{2}nx^2} dx 
		 	\geq
		 	\frac{n^{1/2}} {\sqrt{2 \pi}} \int\limits_{U_n(\frac{1}{2}n^{-1/2+\varepsilon})} e^{-\frac{1}{2}\,\vec{\theta}^T \hat{Q}\vec{\theta}} d \vec{\theta}.
		\end{aligned}
	\end{equation}
	Combining  (\ref{Integral_hatQ}), (\ref{1/2U_n}) and (\ref{LcapU_n->Un}) we obtain (\ref{eq_Lemma_exp-V0}).
\end{Proof}

To prove  Lemma \ref{Lemma_Integral} we separate the integrand in (\ref{Integral_Lemma_Integral}) into three factors. 
\begin{itemize}
	\item $\displaystyle \exp\left(i\,\vec{\theta}^T Q \vec{\alpha_n}\right)$ --- the oscillatory factor, 
	\item  $\displaystyle \exp\left(a\sum\limits_{(v_j,v_k)\in EG}\Delta_{jk}^2\right)$ --- the regular factor,
	\item $\displaystyle \exp\left(b\sum\limits_{(v_j,v_k)\in EG} \Delta_{jk}^4 + R_n(\vec{\theta})\right)$ ---  the residual factor.
\end{itemize}
The proof consists of the following steps.
\begin{enumerate}
	\item In Lemma \ref{Lemma_x^4_Integral} we estimate an integral analogous to (\ref{Integral_Lemma_Integral}) 
	but without an oscillatory factor.
	\item Using Lemma \ref{Lemma_x_Integral}, we get rid of the oscillatory factor in (\ref{Integral_Lemma_Integral}). 
	\item Combining Lemma \ref{Lemma_x_Integral} and Lemma \ref{Lemma_x^4_Integral}, we complete the proof of Lemma \ref{Lemma_Integral}.
\end{enumerate}
At first, we prove two technical statements.
%		J_n = \int\limits_{L \cap V_0}	
%			\exp\left(
%			i\,\vec{\theta}^T Q \vec{\alpha_n} - a\sum\limits_{(v_j,v_k)\in EG}\Delta_{ij}^2
%			- b\sum\limits_{(v_j,v_k)\in EG} \Delta_{jk}^4 + R_n(\vec{\theta})
%			\right)
%		 dL.

\begin{Lemma}\label{Lemma_x^4_Int}
	For any $a>0$ and sequence of functions $r_n(x)$ such that as $n\rightarrow \infty$
	\begin{equation}\label{Lemma_x^4_Int_condit}
		\sup_{|x|\leq n^{-1/2+\varepsilon}}|r_n(x)| = o(1).
	\end{equation}
	Then as $n\rightarrow \infty$
	\begin{equation}\label{eq_Lemma_x^2_Int}
		\int\limits_{-n^{-1/2+\varepsilon}}\limits^{n^{-1/2+\varepsilon}}nx^2 e^{-anx^2 + r_n(x)} dx = \left( \frac{1}{2a} +o(1)\right)
		\int\limits_{-n^{-1/2+\varepsilon}}\limits^{n^{-1/2+\varepsilon}} e^{-anx^2 + r_n(x)} dx.
	\end{equation}
	and
		\begin{equation}\label{eq_Lemma_x^4_Int}
		\int\limits_{-n^{-1/2+\varepsilon}}\limits^{n^{-1/2+\varepsilon}}n^2x^4 e^{-anx^2 + r_n(x)} dx = \left( \frac{3}{4a^2} +o(1)\right)
		\int\limits_{-n^{-1/2+\varepsilon}}\limits^{n^{-1/2+\varepsilon}} e^{-anx^2 + r_n(x)} dx.
	\end{equation}
\end{Lemma}
\begin{Proof} {\it Lemma \ref{Lemma_x^4_Int}.}
	Note that
	\begin{equation}\label{Lemma_x^2_Int_1}
		\int\limits_{-n^{\varepsilon}}\limits^{n^{\varepsilon}} t^2 e^{-at^2} dt
		= \left(1 + O\left(\exp(-cn^{2\varepsilon})\right)\right)
		\int\limits_{-\infty}\limits^{+\infty} t^2 e^{-at^2} dt 
		= \left(\frac{1}{2a} + o(1)\right)\int\limits_{-n^{\varepsilon}}\limits^{n^{\varepsilon}} e^{-at^2} dt
	\end{equation}
	and
	\begin{equation}\label{Lemma_x^4_Int_1}
		\int\limits_{-n^{\varepsilon}}\limits^{n^{\varepsilon}} t^4 e^{-at^2} dt
		= \left(1 + O\left(\exp(-cn^{2\varepsilon})\right)\right)
		\int\limits_{-\infty}\limits^{+\infty} t^4 e^{-at^2} dt 
		= \left(\frac{3}{4a^2} + o(1)\right)\int\limits_{-n^{\varepsilon}}\limits^{n^{\varepsilon}} e^{-at^2} dt.
	\end{equation}
	Using (\ref{Lemma_x^4_Int_condit}), we get that as $n\rightarrow \infty$
	\begin{equation}\label{Lemma_x^4_Int_2}  
		\sup_{|x|\leq n^{-1/2+\varepsilon}}|e^{r_n(x)} - 1| = o(1).
	\end{equation}
	 Making the substitution $t = \sqrt{n}x$ and combining (\ref{Lemma_x^2_Int_1}) and (\ref{Lemma_x^4_Int_1}) with (\ref{Lemma_x^4_Int_2}), 
	 we obtain (\ref{eq_Lemma_x^2_Int}) and (\ref{eq_Lemma_x^4_Int}), respectively. 
\end{Proof}

\begin{Lemma}\label{Lemma_x_Int}
	Under the assumptions of Lemma \ref{Lemma_x^4_Int}, let $r_n(x)$ be differentiable and as $n\rightarrow \infty$
	\begin{equation}\label{Lemma_x_Int_condit}
		\sup_{|x|\leq n^{-1/2+\varepsilon}}|r'_n(x)| = O(n^{-1/2+3\varepsilon}).
	\end{equation}
	Then as $n\rightarrow \infty$
	\begin{equation}\label{eq_Lemma_x_Int}
		\int\limits_{-n^{-1/2+\varepsilon}}\limits^{n^{-1/2+\varepsilon}} x e^{-anx^2 + r_n(x)} dx = O\left(n^{-3/2+4\varepsilon}\right)
		\int\limits_{-n^{-1/2+\varepsilon}}\limits^{n^{-1/2+\varepsilon}} e^{-anx^2 + r_n(x)} dx .
	\end{equation}
\end{Lemma}
\begin{Proof} {\it Lemma \ref{Lemma_x_Int}.}
	Note that
	\begin{equation}\label{Lemma_x_Int_1}
		\int\limits_{0}\limits^{n^{\varepsilon}} t e^{-at^2} dt
		= \left(1 + O\left(\exp(-cn^{2\varepsilon})\right)\right)
		\int\limits_{0}\limits^{+\infty} t e^{-at^2} dt 
		= \frac{1}{2a} + O\left(\exp(-cn^{2\varepsilon})\right).
	\end{equation}
	According to the Mean Value Theorem, we have that for some $|\tilde{x}|\leq |x|$ 
	\begin{equation}
		|e^{r_n(x)} - e^{r_n(-x)}|  
		= |e^{r_n(\tilde{x})} r'_n(\tilde{x}) 2x |.
	\end{equation}
	Using (\ref{Lemma_x_Int_condit}),  we get that as $n\rightarrow \infty$
	\begin{equation}\label{Lemma_x_Int_2}  
		\sup_{|x|\leq n^{-1/2+\varepsilon}} |e^{r_n(x)} - e^{r_n(-x)}| = O\left(n^{-1+4\varepsilon}\right).
	\end{equation}
	We have that
	\begin{equation}
		\int\limits_{-n^{-1/2+\varepsilon}}\limits^{n^{-1/2+\varepsilon}} x e^{-anx^2 + r_n(x)} dx = 
		\int\limits_{0}\limits^{n^{-1/2+\varepsilon}} x e^{-anx^2} \left(e^{r_n(x)} - e^{r_n(-x)}\right)dx.
	\end{equation}
	Making the substitution $t = \sqrt{n}x$ and combining (\ref{Lemma_x_Int_1}) with (\ref{Lemma_x_Int_2}), we obtain (\ref{eq_Lemma_x_Int}).
\end{Proof}

We use notation
\begin{equation}
	\mu_m = \sum\limits_{j=1}\limits^{n} |\phi_j|^m.
\end{equation}
According to the Generalized Mean Inequality, we have that
\begin{equation}\label{mu_1<mu_4}
	\mu_1/n \leq 
%	\left(\mu_2/n\right)^{1/2} \leq 
\left(\mu_4/n\right)^{1/4}.
\end{equation}
Since
\begin{equation}\label{phi_k_theta_k}
	\phi_k = (d_k+1)\theta_k + \sum\limits_{(v_k,v_j)\notin EG} \theta_j,
\end{equation} 
and (see (\ref{Lemma_minor_d1}))
\begin{equation}
	d_k \geq \sigma n/2	
\end{equation}
we obtain that
\begin{equation}\label{|theta|<|phi|}
	|\theta_k| \leq \frac{2}{\sigma n}\left(|\phi_k| +  \sum\limits_{j\neq k} |\theta_j|\right)
\end{equation}
Using Lemma \ref{Lemma_normQ}, we find that 
\begin{equation}\label{Using_lemmanormQ}
	\sum\limits_{j\neq k} |\theta_j| \leq ||\theta||_1 \leq \frac{c_\infty}{n} ||\phi||_1 = \frac{c_\infty}{n} \mu_1
\end{equation}
 Combining (\ref{mu_1<mu_4}), (\ref{|theta|<|phi|}) and (\ref{Using_lemmanormQ}), we get that
%\begin{equation}
%	\theta_k^2 \leq \frac{4}{\sigma^2 n^2}\left(|\phi_k| +  \frac{c_\infty}{n}\mu_1 \right)^2
%	\leq \frac{4}{\sigma^2 n^2}\left(|\phi_k| +  c_\infty \left(\mu_2/n\right)^{1/2} \right)^2.
%\end{equation}
%and
\begin{equation}
	\theta_k^4 \leq \frac{16}{\sigma^4 n^4}\left(|\phi_k| +  \frac{c_\infty}{n}\mu_1 \right)^4
	\leq \frac{16}{\sigma^4 n^4}\left(|\phi_k| +  c_\infty \left(\mu_4/n\right)^{1/4} \right)^4.
\end{equation}
Using  the inequality $(x+y)^4 \leq 8(x^4 + y^4)$, we obtain that
%\begin{equation}\label{|theta|^2<|phi|^2}
%	\theta_k^2 \leq c_\phi^{(2)} \frac{\phi_k^2}{n^2} + c_\mu^{(2)} \frac{\mu_2}{n^3},
%\end{equation}
%and
\begin{equation}\label{|theta|^4<|phi|^4}
	\theta_k^4 \leq c_\phi \frac{\phi_k^4}{n^4} + c_\mu \frac{\mu_4}{n^5},
\end{equation}
where constants $c_\phi, c_\mu >0$ depend only on $\sigma$.
\begin{Lemma}\label{Lemma_x^4_Integral}
	Let assumptions of Lemma \ref{Lemma_normQ} hold. Let $\{a_n\}$ be  sequence of positive numbers having limit $a>0$. 
	Then for any $b>0$ as $n\rightarrow \infty$
	\begin{equation}\label{eq_Lemma_x^4_Integral}
		\int_{U_n(n^{-1/2+\varepsilon})} \exp\left( - a_n \vec{\theta}^T \hat{Q} \vec{\theta} - b \sum\limits_{(v_j,v_k)\in EG} 
		\Delta_{jk}^4 \right) d\vec{\theta}
		= \Theta_{k_1,k_2} \left(\int_{\mathbb{R}^n} e^{ - a_n \vec{\theta}^T \hat{Q} \vec{\theta} } d\vec{\theta} \right),
	\end{equation}
	where constants $k_1,k_2>0$ depend only on $a$, $b$ and $\sigma$.
\end{Lemma}
\begin{Proof} {\it Lemma \ref{Lemma_x^4_Integral}.}
	Using  the inequality $(x+y)^4 \leq 8(x^4 + y^4)$, we find that
	\begin{equation}
		\sum\limits_{(v_j,v_k)\in EG} 
		\Delta_{jk}^4  \leq 8n \sum\limits_{j=1}\limits^{n} \theta_j^4.
	\end{equation}
	We define $R_1(\vec{\theta}) = 8n \sum\limits_{j=1}\limits^{n} \theta_j^4$. Thus we have that
	\begin{equation}
			\int_{U_n(n^{-1/2+\varepsilon})} \exp\left( - a_n \vec{\theta}^T \hat{Q} \vec{\theta} - b \sum\limits_{(v_j,v_k)\in EG} 
		\Delta_{jk}^4 \right) d\vec{\theta} \geq \int_{\mathbb{R}^n} e^{ - a_n \vec{\theta}^T \hat{Q} \vec{\theta} - R_1(\vec{\theta})} d\vec{\theta}.
	\end{equation}
	Using (\ref{1/2U_n}), we find that as $n\rightarrow \infty$ 
	\begin{equation}\label{Lemma_x^4_Integral_-1}
		\begin{aligned}
		\int_{U_n(n^{-1/2+\varepsilon})} \phi_1^4 e^{ - a_n \vec{\theta}^T \hat{Q} \vec{\theta} - R_1(\vec{\theta})} d\vec{\theta} 
		= \int_{U_n(\frac{4}{\sigma}n^{-1/2+\varepsilon})} \phi_1^4 e^{ - a_n \vec{\theta}^T \hat{Q} \vec{\theta} - R_1(\vec{\theta})} d\vec{\theta}+\\
		+ O\left(\exp(-cn^{2\varepsilon})\right)\int_{\mathbb{R}^n} e^{ - a_n \vec{\theta}^T \hat{Q} \vec{\theta} } d\vec{\theta}
		\end{aligned}
	\end{equation}
	for some $c>0$ depending only on $a$ and $\sigma$. It follows that
	\begin{equation}\label{Lemma_x^4_Integral_0}
		\begin{aligned}	
				\int\limits_{-n^{-1/2+\varepsilon}}\limits^{n^{-1/2+\varepsilon}}\cdots \int\limits_{-n^{-1/2+\varepsilon}}\limits^{n^{-1/2+\varepsilon}}
				\left( \int \limits_{-\frac{4}{\sigma}n^{-1/2+\varepsilon}}\limits^{\frac{4}{\sigma}n^{-1/2+\varepsilon}} 
				\phi_1^4 e^{ - a_n \vec{\theta}^T \hat{Q} \vec{\theta} - R_1(\vec{\theta})} d \theta_1  \right)
				d \theta_2 \ldots d \theta_n = \\
				= \int_{U_n(n^{-1/2+\varepsilon})} \phi_1^4 e^{ - a_n \vec{\theta}^T \hat{Q} \vec{\theta} - R_1(\vec{\theta})} d\vec{\theta} 
				+ O\left(\exp(-cn^{2\varepsilon})\right)\int_{\mathbb{R}^n} e^{ - a_n \vec{\theta}^T \hat{Q} \vec{\theta} } d\vec{\theta}
		\end{aligned}
	\end{equation}
	Using (\ref{g_theta}), we find that
	\begin{equation}\label{Lemma_x^4_Integral_1}
		\begin{aligned}	
				\int\limits_{-n^{-1/2+\varepsilon}}\limits^{n^{-1/2+\varepsilon}}\cdots \int\limits_{-n^{-1/2+\varepsilon}}\limits^{n^{-1/2+\varepsilon}}
				\left( \int \limits_{-\frac{4}{\sigma}n^{-1/2+\varepsilon}}\limits^{\frac{4}{\sigma}n^{-1/2+\varepsilon}} 
				\phi_1^4 e^{ - a_n \vec{\theta}^T \hat{Q} \vec{\theta} - R_1(\vec{\theta})} d \theta_1  \right)
				d \theta_2 \ldots d \theta_n = \\
				= \int\limits_{-n^{-1/2+\varepsilon}}\limits^{n^{-1/2+\varepsilon}}\cdots \int\limits_{-n^{-1/2+\varepsilon}}\limits^{n^{-1/2+\varepsilon}}
				e^{ - a_n g_1(\theta_2,\ldots,\theta_n) - R_2(\vec{\theta})}
				\left( \int \limits_{-\frac{4}{\sigma}n^{-1/2+\varepsilon}}\limits^{\frac{4}{\sigma}n^{-1/2+\varepsilon}} 
				\phi_1^4 e^{ - a_n \frac{\phi_1^2}{d_1+1} - 8bn\theta_1^4} d \theta_1  \right)
				d \theta_2 \ldots d \theta_n,
		\end{aligned}
	\end{equation} 
	where $R_2(\vec{\theta}) = 8n \sum\limits_{j=2}\limits^{n} \theta_j^4$.
	Using (\ref{|theta|<|phi|}), we get that as $n\rightarrow \infty$
	\begin{equation}\label{Lemma_x^4_Integral_2}
		\int \limits_{-\frac{4}{\sigma}n^{-1/2+\varepsilon}}\limits^{\frac{4}{\sigma}n^{-1/2+\varepsilon}} 
				\phi_1^4 e^{ - a_n \frac{\phi_1^2}{d_1+1} - 8bn\theta_1^4} d \theta_1 =  
				\left( 1+ O\left(\exp(-cn^{2\varepsilon})\right) \right)
		\int \limits_{|\phi_1|\leq n^{1/2 + \varepsilon}}  
				\phi_1^4 e^{ - a_n \frac{\phi_1^2}{d_1+1} - 8bn\theta_1^4} d \theta_1
	\end{equation}
	Combining (\ref{Lemma_x^4_Integral_0}), (\ref{Lemma_x^4_Integral_1}), (\ref{Lemma_x^4_Integral_2}) and Lemma \ref{Lemma_x^4_Int} 
	with $x = \phi_1/n$, 
	we obtain that as $n\rightarrow \infty$
	\begin{equation}\label{Lemma_x^4_Integral_3}
		\begin{aligned}
		\int_{U_n(n^{-1/2+\varepsilon})} \phi_1^4 e^{ - a_n \vec{\theta}^T \hat{Q} \vec{\theta} - R_1(\vec{\theta})} d\vec{\theta} 
			\leq c' n^2  \int_{U_n(n^{-1/2+\varepsilon})} e^{ - a_n \vec{\theta}^T \hat{Q} \vec{\theta} - R_1(\vec{\theta})} d\vec{\theta} +\\+ 
			O\left(\exp(-cn^{2\varepsilon})\right)\int_{\mathbb{R}^n} e^{ - a_n \vec{\theta}^T \hat{Q} \vec{\theta} } d\vec{\theta}
		\end{aligned}
	\end{equation} 
	for some constants $c,c'>0$ depending only on $a$ and $\sigma$.

Combining similar to (\ref{Lemma_x^4_Integral_3}) inequalities for $\phi_1, \phi_2, \ldots ,\phi_n$ and using (\ref{|theta|^4<|phi|^4}), 
we find that as $n\rightarrow \infty$ 
	\begin{equation}\label{Lemma_x^4_Integral_4}
		\begin{aligned}
		\int_{U_n(n^{-1/2+\varepsilon})} \theta_1^4 e^{ - a_n \vec{\theta}^T \hat{Q} \vec{\theta} - R_1(\vec{\theta})} d\vec{\theta} 
			\leq \frac{(c_\phi+c_\mu)c'}{n^2}  \int_{U_n(n^{-1/2+\varepsilon})} e^{ - a_n \vec{\theta}^T \hat{Q} \vec{\theta} - R_1(\vec{\theta})} d\vec{\theta} +\\+ 
			O\left(\exp(-cn^{2\varepsilon})\right)\int_{\mathbb{R}^n} e^{ - a_n \vec{\theta}^T \hat{Q} \vec{\theta} } d\vec{\theta}
		\end{aligned}
	\end{equation} 
	for some $c>0$ depending only on $a$ and $\sigma$. Note that as $n\rightarrow \infty$
	\begin{equation}\label{Lemma_x^4_Integral_5}
		\begin{aligned}
			\int\limits_{-n^{-1/2+\varepsilon}}\limits^{n^{-1/2+\varepsilon}}\cdots \int\limits_{-n^{-1/2+\varepsilon}}\limits^{n^{-1/2+\varepsilon}}
				e^{ - a_n g_1(\theta_2,\ldots,\theta_n) - R_2(\vec{\theta})}
				\left( \int \limits_{-n^{-1/2+\varepsilon}}\limits^{n^{-1/2+\varepsilon}} 
				 e^{ - a_n \frac{\phi_1^2}{d_1+1} - 8bn\theta_1^4} d \theta_1  \right)
				d \theta_2 \ldots d \theta_n = \\
				= 
			%	\int\limits_{-n^{-1/2+\varepsilon}}\limits^{n^{-1/2+\varepsilon}}
				\cdots 
		%		\int\limits_{-n^{-1/2+\varepsilon}}\limits^{n^{-1/2+\varepsilon}}
			%	e^{ - a_n \vec{\theta}^T \hat{Q} \vec{\theta} - R_2(\vec{\theta})}
				\left( \int \limits_{-n^{-1/2+\varepsilon}}\limits^{n^{-1/2+\varepsilon}} 
				 e^{ - a_n \frac{\phi_1^2}{d_1+1}} \left(1 - 8bn\theta_1^4 + O\left(n^{-2+8\varepsilon}\right)\right) d \theta_1  \right) 
				d \theta_2 \ldots d \theta_n		
		\end{aligned}
	\end{equation}
	Combining (\ref{Lemma_x^4_Integral_4}) and (\ref{Lemma_x^4_Integral_5}), we get that as $n\rightarrow \infty$ 
	\begin{equation}\label{Lemma_x^4_Integral_6}
		\begin{aligned}
		\int_{U_n(n^{-1/2+\varepsilon})}  e^{ - a_n \vec{\theta}^T \hat{Q} \vec{\theta} - R_1(\vec{\theta})} d\vec{\theta} 
			\geq \left(1 +  \frac{\tilde{c}}{n} \right)  \int_{U_n(n^{-1/2+\varepsilon})} e^{ - a_n \vec{\theta}^T \hat{Q} \vec{\theta} - R_2(\vec{\theta})} d\vec{\theta} +\\+ 
			O\left(\exp(-cn^{2\varepsilon})\right)\int_{\mathbb{R}^n} e^{ - a_n \vec{\theta}^T \hat{Q} \vec{\theta} } d\vec{\theta},
		\end{aligned}
	\end{equation} 
	where $\tilde{c}$ depends only on $a$, $b$ and $\sigma$.
	
		We continue similarly to (\ref{Lemma_x^4_Integral_6})  
	\begin{equation}\label{Lemma_x^4_Integral_7}
		\begin{aligned}
		\int_{U_n(n^{-1/2+\varepsilon})}  e^{ - a_n \vec{\theta}^T \hat{Q} \vec{\theta} - R_k(\vec{\theta})} d\vec{\theta} 
			\geq \left(1 +  \frac{\tilde{c}}{n} \right)  \int_{U_n(n^{-1/2+\varepsilon})} e^{ - a_n \vec{\theta}^T \hat{Q} \vec{\theta} - R_{k+1}(\vec{\theta})} 		     
			d\vec{\theta} 			+\\+ 
			O\left(\exp(-cn^{2\varepsilon})\right)\int_{\mathbb{R}^n} e^{ - a_n \vec{\theta}^T \hat{Q} \vec{\theta} } d\vec{\theta},
		\end{aligned}
	\end{equation} 
		where 
		\begin{equation}
			R_k(\vec{\theta}) = 8n \sum\limits_{j=k}\limits^{n} \theta_j^4.
		\end{equation}
		Combining all inequalities of (\ref{Lemma_x^4_Integral_7}) for $R_1, R_2, \ldots , R_n$, we get that
		\begin{equation}\label{Lemma_x^4_Integral_8}
			\begin{aligned}
				\int_{U_n(n^{-1/2+\varepsilon})}  e^{ - a_n \vec{\theta}^T \hat{Q} \vec{\theta} - R_1(\vec{\theta})} d\vec{\theta} 
				\geq \left(1 +  \frac{\tilde{c}}{n} \right)^n 
				\int_{U_n(n^{-1/2+\varepsilon})}  e^{ - a_n \vec{\theta}^T \hat{Q} \vec{\theta}} d\vec{\theta} +\\
				+O\left(\exp(-cn^{2\varepsilon})\right)\int_{\mathbb{R}^n} e^{ - a_n \vec{\theta}^T \hat{Q} \vec{\theta} } d\vec{\theta}
			\end{aligned}
		\end{equation}
		for some $c>0$ depending only on $a$ and $\sigma$. Note also that
		\begin{equation}\label{Lemma_x^4_Integral_9}
			\int_{U_n(n^{-1/2+\varepsilon})}  e^{ - a_n \vec{\theta}^T \hat{Q} \vec{\theta} - R_1(\vec{\theta})} d\vec{\theta} 
			\leq \int_{U_n(n^{-1/2+\varepsilon})}  e^{ - a_n \vec{\theta}^T \hat{Q} \vec{\theta}} d\vec{\theta}
		\end{equation}
	Combining (\ref{1/2U_n}), (\ref{Lemma_x^4_Integral_8}) and (\ref{Lemma_x^4_Integral_9}), we obtain (\ref{eq_Lemma_x^4_Integral}).
	\end{Proof}

\begin{Lemma}\label{Lemma_phi_Integral}
	Let the assumptions of Lemma \ref{Lemma_normQ} hold. 
	For positive constants $a,b,d_1,d_2$ let 
	sequence of differentiable functions $\{R_n(\vec{\theta})\}$ be such that 
	\begin{equation}
		Re\left(R_n(\vec{\theta})\right) \leq d_1\frac{\vec{\theta}^T \hat{Q}\vec{\theta}}{n},
	\end{equation}
	and for $\vec{\theta} \in U_n(\frac{4}{\sigma}n^{-1/2+\varepsilon})$
	\begin{equation}
	 \left|\frac{\partial R_n(\vec{\theta})}{\partial\theta_k}\right| \leq d_2 \,n^{-1/2+3\varepsilon}.
	\end{equation}
	Then as $n\rightarrow\infty$ 
	\begin{equation}\label{eq_Lemma_phi_Integral1}
		\begin{aligned}
		\int\limits_{U_n(n^{-1/2+\varepsilon})}	
			\phi_k
			e^{
			i \frac{b}{n}\phi_k - a\, \vec{\theta}^T \hat{Q} \vec{\theta} + R_n(\vec{\theta})
			}
		 d\vec{\theta}  = \Theta_{k_1,k_2}\left(
		 \int\limits_{U_n(n^{-3/2+\varepsilon})}	
			e^{
			  i \frac{b}{n}\phi_k - a\, \vec{\theta}^T \hat{Q} \vec{\theta} + R_n(\vec{\theta})
			}
		 d\vec{\theta}\right) +\\+ O(n^{-1/2+4\varepsilon}) 
		 \int\limits_{\mathbb{R}^n}	
			e^{
			 - a \, \vec{\theta}^T \hat{Q} \vec{\theta}}
			 d\vec{\theta}
			 \end{aligned}
	\end{equation}
		and
		\begin{equation}\label{eq_Lemma_phi_Integral2}
		\begin{aligned}
		\int\limits_{U_n(n^{-1/2+\varepsilon})}	
			\phi_k^2
			e^{
			i \frac{b}{n}\phi_k - a\, \vec{\theta}^T \hat{Q} \vec{\theta} + R_n(\vec{\theta})
			}
		 d\vec{\theta}  = \Theta_{k_1,k_2}\left( n
		 \int\limits_{U_n(n^{-3/2+\varepsilon})}	
			e^{
			  i \frac{b}{n}\phi_k - a\, \vec{\theta}^T \hat{Q} \vec{\theta} + R_n(\vec{\theta})
			}
		 d\vec{\theta}\right) +\\+ O\left(\exp(-\tilde{c}n^{2\varepsilon})\right) 
		 \int\limits_{\mathbb{R}^n}	
			e^{
			 - a \, \vec{\theta}^T \hat{Q} \vec{\theta}}
			 d\vec{\theta},
			 \end{aligned}
	\end{equation}	
	where constants $k_1,k_2, \tilde{c} >0$ depend only on $a$, $b$, $d_1$, $d_2$ and $\sigma$ .
\end{Lemma}

\begin{Proof} {\it Lemma \ref{Lemma_phi_Integral}.}
For our purpose it is convenient to assume that $k = 1$. Note that
	\begin{equation}
		\int\limits_{U_n(n^{-1/2+\varepsilon})}	
			|e^{
			i \frac{b}{n}\phi_1 - a\, \vec{\theta}^T \hat{Q} \vec{\theta} + R_n(\vec{\theta})
			}|
			d\vec{\theta} \leq \int\limits_{U_n(n^{-1/2+\varepsilon})}	
			e^{
			- a\, \vec{\theta}^T \hat{Q} \vec{\theta} +  d_1\frac{\vec{\theta}^T \hat{Q}\vec{\theta}}{n}
			}
			d\vec{\theta}.
	\end{equation} 
	Using (\ref{Integral_hatQ}), we get that as $n\rightarrow \infty$
	\begin{equation}\label{Lemma_phi_Integral_-1}
		\int\limits_{U_n(n^{-1/2+\varepsilon})}	
			|e^{
			i \frac{b}{n}\phi_1 - a\, \vec{\theta}^T \hat{Q} \vec{\theta} + R_n(\vec{\theta})
			}|
			d\vec{\theta} = O\left(1\right) \int\limits_{U_n(n^{-1/2+\varepsilon})}
			e^{
			- a\, \vec{\theta}^T \hat{Q} \vec{\theta}
			}
			d\vec{\theta}
			.
	\end{equation}
	Similar to (\ref{Lemma_x^4_Integral_-1}), for $m = 1,2$ we find that as $n\rightarrow \infty$ 
	\begin{equation}
		\begin{aligned}
		\int_{U_n(\frac{4}{\sigma}n^{-1/2+\varepsilon}) - U_n(n^{-1/2+\varepsilon})} \phi_1^m 
		&|e^{
		i \frac{b}{n}\phi_1 - a\, \vec{\theta}^T \hat{Q} \vec{\theta} + R_n(\vec{\theta})
		}| 
		 d\vec{\theta} 
		=\\&= O\left(\exp(-cn^{2\varepsilon})\right)\int_{\mathbb{R}^n} e^{ - a\, \vec{\theta}^T \hat{Q} \vec{\theta} } d\vec{\theta}
		\end{aligned}
	\end{equation}
	for some $c>0$ depending only on $a,d_1$ and $\sigma$. It follows that
	\begin{equation}\label{Lemma_phi_Integral_0}
		\begin{aligned}	
				\int\limits_{-n^{-1/2+\varepsilon}}\limits^{n^{-1/2+\varepsilon}}\cdots \int\limits_{-n^{-1/2+\varepsilon}}\limits^{n^{-1/2+\varepsilon}}
				\left( \int \limits_{-\frac{4}{\sigma}n^{-1/2+\varepsilon}}\limits^{\frac{4}{\sigma}n^{-1/2+\varepsilon}} 
				\phi_1^m 
				e^{
				i \frac{b}{n}\phi_1 - a\, \vec{\theta}^T \hat{Q} \vec{\theta} + R_n(\vec{\theta})
				} d \theta_1  \right)
				d \theta_2 \ldots d \theta_n = \\
				= \int_{U_n(n^{-1/2+\varepsilon})} \phi_1^m 
				e^{
				i \frac{b}{n}\phi_1 - a\, \vec{\theta}^T \hat{Q} \vec{\theta} + R_n(\vec{\theta})
				} 
				d\vec{\theta} 
				+ O\left(\exp(-cn^{2\varepsilon})\right)\int_{\mathbb{R}^n} e^{ - a\, \vec{\theta}^T \hat{Q} \vec{\theta} } d\vec{\theta}
		\end{aligned}
	\end{equation}
	We define 
	\begin{equation}
			R_n'(\vec{\theta}) = R_n(0,\theta_2, \ldots,\theta_{n})
	\end{equation}
	and 
	\begin{equation}
			r_1(\vec{\theta}) = R_n(\vec{\theta}) - R_n'(\vec{\theta}). 
	\end{equation} 
	According to Mean Value Theorem, for $\theta \in U_n(\frac{4}{\sigma}n^{-1/2+\varepsilon})$ we have that 
	\begin{equation}
		|r_1(\vec{\theta})| = |R_n(\vec{\theta}) - R_n'(\vec{\theta})|
		= \left|\theta_1 \frac{\partial R_n(\tilde{\vec{\theta}})}{\partial\theta_k}\right| = O(n^{-1+4\varepsilon}).
	\end{equation}
	Using (\ref{g_theta}), we find that
	\begin{equation}\label{Lemma_phi_Integral_1}
		\begin{aligned}	
				\int\limits_{-n^{-1/2+\varepsilon}}\limits^{n^{-1/2+\varepsilon}}
				\cdots \int\limits_{-n^{-1/2+\varepsilon}}\limits^{n^{-1/2+\varepsilon}}
				\left( \int \limits_{-\frac{4}{\sigma}n^{-1/2+\varepsilon}}\limits^{\frac{4}{\sigma}n^{-1/2+\varepsilon}} 
				\phi_1^m 
				e^{
				i \frac{b}{n}\phi_1 - a\, \vec{\theta}^T \hat{Q} \vec{\theta} + R_n(\vec{\theta})
				} d \theta_1  \right)
				d \theta_2 \ldots d \theta_n = \\
				= \cdots \int\limits_{-n^{-1/2+\varepsilon}}\limits^{n^{-1/2+\varepsilon}}
				e^{ - a g_1(\theta_2,\ldots,\theta_n) + R_n'(\vec{\theta})}
				\left( \int \limits_{-\frac{4}{\sigma}n^{-1/2+\varepsilon}}\limits^{\frac{4}{\sigma}n^{-1/2+\varepsilon}} 
				\phi_1^m 
				e^{i \frac{b}{n}\phi_1 - a\, \frac{\phi_1^2}{d_1+1} + r_1(\vec{\theta})} 
				d \theta_1  \right) d \theta_2 \ldots d \theta_n,
		\end{aligned}
	\end{equation} 
	Using (\ref{|theta|<|phi|}), we get that as $n\rightarrow \infty$
	\begin{equation}\label{Lemma_phi_Integral_2}
		\begin{aligned}
		\int \limits_{-\frac{4}{\sigma}n^{-1/2+\varepsilon}}\limits^{\frac{4}{\sigma}n^{-1/2+\varepsilon}}& 
				\phi_1^m e^{i \frac{b}{n}\phi_1 - a\, \frac{\phi_1^2}{d_1+1} + r_1(\vec{\theta})} d \theta_1 =  \\ &=
				\left( 1+ O\left(\exp(-cn^{2\varepsilon})\right) \right)
		\int \limits_{|\phi_1|\leq n^{1/2 + \varepsilon}}  
				\phi_1^m e^{i \frac{b}{n}\phi_1 - a\, \frac{\phi_1^2}{d_1+1} + r_1(\vec{\theta})} d \theta_1
		\end{aligned}
	\end{equation}
	Combining (\ref{Lemma_phi_Integral_0}), (\ref{Lemma_phi_Integral_1}), (\ref{Lemma_phi_Integral_2}) with $m = 2$ and (\ref{eq_Lemma_x^2_Int})
	with $x = \phi_1/n$, 
	we obtain (\ref{eq_Lemma_phi_Integral2}).
	
	Note that
	\begin{equation}\label{Lemma_phi_Integral_3}
		\begin{aligned}
		\int \limits_{|\phi_1|\leq n^{1/2 + \varepsilon}}  
				\phi_1 &e^{i \frac{b}{n}\phi_1 - a\, \frac{\phi_1^2}{d_1+1} + r_1(\vec{\theta})} d \theta_1
				=\\ &=\int \limits_{|\phi_1|\leq n^{1/2 + \varepsilon}}  
				\phi_1 \left(1 + i \frac{b}{n}\phi_1 + O\left(n^{-1+2\varepsilon}\right)\right)e^{ - a\, \frac{\phi_1^2}{d_1+1} + r_1(\vec{\theta})} d \theta_1
		\end{aligned}
	\end{equation}
	Since $\partial r_1/ \partial \theta_1 = \partial R_n/ \partial \theta_1$, using (\ref{eq_Lemma_x_Int})
	with $x = \phi_1/n$, we get that
	\begin{equation}\label{Lemma_phi_Integral_4}
		\begin{aligned}
				\int \limits_{|\phi_1|\leq n^{1/2 + \varepsilon}}  
				\phi_1 e^{ - a\, \frac{\phi_1^2}{d_1+1} + r_1(\vec{\theta})} d \theta_1 
				= O\left(n^{-1/2+4\varepsilon}\right) 			
				\int \limits_{|\phi_1|\leq n^{1/2 + \varepsilon}}  
				 e^{ - a\, \frac{\phi_1^2}{d_1+1} + r_1(\vec{\theta})} d \theta_1 
				 %=\\
%				 &= O\left(n^{-1/2+4\varepsilon}\right) 			
%				\int \limits_{|\phi_1|\leq n^{1/2 + \varepsilon}}  
%				 e^{i \frac{b}{n}\phi_1 - a\, \frac{\phi_1^2}{d_1+1} + r_1(\vec{\theta})} d \theta_1.
		\end{aligned}
	\end{equation}
	
	Combining (\ref{Lemma_phi_Integral_0}), (\ref{Lemma_phi_Integral_1}), (\ref{Lemma_phi_Integral_2}) with $m = 1$ and 
	(\ref{eq_Lemma_phi_Integral2}), (\ref{Lemma_phi_Integral_-1}) with $b = 0$ and (\ref{Lemma_phi_Integral_4}), we
	obtain (\ref{eq_Lemma_phi_Integral1}). 
\end{Proof}

\begin{Lemma}\label{Lemma_x_Integral}
	Let the assumptions of Lemma \ref{Lemma_normQ} hold. 
	For positive constants $a,b,d_1,d_2$ let sequence of vectors $\{\vec{\beta}_n\}$ and 
	sequence of real differentiable functions $\{R_n(\vec{\theta})\}$ be such that 
	\begin{equation}
		||\vec{\beta}_n||_{\infty} \leq b,
	\end{equation}
	\begin{equation}
		R_n(\vec{\theta}) \leq d_1\frac{\vec{\theta}^T \hat{Q}\vec{\theta}}{n},
	\end{equation}
	and for $\vec{\theta} \in U_n(\frac{4}{\sigma}n^{-1/2+\varepsilon})$
	\begin{equation}
	 \left\|\frac{\partial R_n(\vec{\theta})}{\partial\vec{\theta}}\right\|_{\infty} \leq d_2 \,n^{-1/2+3\varepsilon}.
	\end{equation}
	Then as $n\rightarrow\infty$ 
	\begin{equation}\label{eq_Lemma_x_Integral}
		\begin{aligned}
		\int\limits_{U_n(n^{-1/2+\varepsilon})}	
			e^{
			i\, \vec{\beta}_n^T \vec{\theta} - a\, \vec{\theta}^T \hat{Q} \vec{\theta} + R_n(\vec{\theta})
			}
		 d\vec{\theta}  =
		 \Theta_{\tilde{k}_1,\tilde{k}_2} \left(
		 \int\limits_{U_n(n^{-1/2+\varepsilon})}	
			e^{
			 - a \, \vec{\theta}^T \hat{Q} \vec{\theta} + \tilde{R}_n(\vec{\theta})
			}
		 d\vec{\theta}
		 \right)
		  +\\+ O(n^{-1/2+4\varepsilon}) 
		 \int\limits_{\mathbb{R}^n}	
			e^{
			 - a \, \vec{\theta}^T \hat{Q} \vec{\theta}}
			 d\vec{\theta},
			 \end{aligned}
	\end{equation}
	where $\tilde{R}_n = R_n - \frac{1}{2}\sum\limits_{j=1}\limits^{n} \beta_j^2\theta_j^2$ 
	and constants $\tilde{k}_1,\tilde{k}_2$ depend only on $a$, $b$, $d_1$, $d_2$ and $\sigma$ .
\end{Lemma}
\begin{Proof} {\it Lemma \ref{Lemma_x_Integral}.}
	Using (\ref{Lemma_phi_Integral_-1}), we get that as $n \rightarrow \infty$
	\begin{equation}\label{x_Integral1}
		\begin{aligned}
			&\int\limits_{U_n(n^{-1/2+\varepsilon})}	
				e^{
				i\, \vec{\beta}_n^T \vec{\theta} - a\, \vec{\theta}^T \hat{Q} \vec{\theta} + R_n(\vec{\theta})
				}
			 d\vec{\theta} 
			 =\\&= 
			 	 \int\limits_{U_n(n^{-1/2+\varepsilon})}	
			 \left(1 + i\beta_1\theta_1 -\frac{\beta_1^2\theta_1^2}{2} + O(n^{-3/2+3\varepsilon}) \right)
						e^{
				i\, (\vec{\beta}_n^T \vec{\theta}- \beta_1\theta_1) - a\, \vec{\theta}^T \hat{Q} \vec{\theta} + R_n(\vec{\theta})
				} d\vec{\theta}	
				=\\ 
				&= \int\limits_{U_n(n^{-1/2+\varepsilon})}	
				e^{
				i\, (\vec{\beta}_n^T \vec{\theta}- \beta_1\theta_1) - a\, \vec{\theta}^T \hat{Q} \vec{\theta} + R_n(\vec{\theta}) - \frac{1}{2}\beta_1^2\theta_1^2
				} d\vec{\theta}	+ O(n^{-3/2+3\varepsilon})
				\int\limits_{U_n(n^{-1/2+\varepsilon})}	
				e^{
				 - a\, \vec{\theta}^T \hat{Q} \vec{\theta} 
				} d\vec{\theta}	+				
				\\ 
				&+  
			 \int\limits_{U_n(n^{-1/2+\varepsilon})}	i\beta_1\theta_1
				e^{
				i\, (\vec{\beta}_n^T \vec{\theta}- \beta_1\theta_1) - a\, \vec{\theta}^T \hat{Q} \vec{\theta} + R_n(\vec{\theta})
				} d\vec{\theta}.	
		\end{aligned}	 
	\end{equation}
	Taking into account (\ref{phi_k_theta_k}) and using Lemma \ref{Lemma_phi_Integral}, we find that as $n\rightarrow \infty$
	\begin{equation}\label{x_Integral2}
		\begin{aligned}
			\int\limits_{U_n(n^{-1/2+\varepsilon})}	
				\phi_k e^{
				i\, (\vec{\beta}_n^T \vec{\theta}- \beta_1\theta_1) - a\, \vec{\theta}^T \hat{Q} \vec{\theta} + R_n(\vec{\theta})
				} &d\vec{\theta} 
				=\\
				=
				\Theta_{k_1,k_2}\left(
		 \int\limits_{U_n(n^{-3/2+\varepsilon})}	
			e^{
			  i (\vec{\beta}_n^T \vec{\theta}- \beta_1\theta_1)  - a\, \vec{\theta}^T \hat{Q} \vec{\theta} + R_n(\vec{\theta})
			}
		 d\vec{\theta}\right)&+\\ + O(n^{-1/2+4\varepsilon})& 
		 \int\limits_{\mathbb{R}^n}	
			e^{
			 - a \, \vec{\theta}^T \hat{Q} \vec{\theta}}
			 d\vec{\theta}=
			 \\
			 = 
			 \Theta_{k_1,k_2}\left(
		 \int\limits_{U_n(n^{-3/2+\varepsilon})}	
			e^{
			  i (\vec{\beta}_n^T \vec{\theta}- \beta_1\theta_1)  - a\, \vec{\theta}^T \hat{Q} \vec{\theta} + R_n(\vec{\theta}) - \frac{1}{2}\beta_1^2\theta_1^2
			}
		 d\vec{\theta}\right)&  +\\+O(n^{-1/2+4\varepsilon})& 
		 \int\limits_{\mathbb{R}^n}	
			e^{
			 - a \, \vec{\theta}^T \hat{Q} \vec{\theta}}
			 d\vec{\theta},
		\end{aligned}
	\end{equation} 
	where constants $k_1,k_2$ depend only on $a$, $b$, $d_1$, $d_2$ and $\sigma$. 
	
	According to Lemma \ref{Lemma_normQ}, we have that
	\begin{equation}
		||\hat{Q}^{-1}||_1 = ||\hat{Q}^{-1}||_\infty \leq \frac{c_\infty}{n}. 
	\end{equation}
	Thus as $n\rightarrow \infty$
	\begin{equation}\label{x_Integral3}
		\begin{aligned}
			&\left|\int\limits_{U_n(n^{-1/2+\varepsilon})}	
				\theta_k e^{
				i\, (\vec{\beta}_n^T \vec{\theta}- \beta_1\theta_1) - a\, \vec{\theta}^T \hat{Q} \vec{\theta} + R_n(\vec{\theta})
				} d\vec{\theta}\right| 
				\leq\\
				&\leq \frac{\tilde{c}}{n}
				\left|\int\limits_{U_n(n^{-3/2+\varepsilon})}	
					e^{
			 		 i(\vec{\beta}_n^T \vec{\theta}- \beta_1\theta_1)  - a\, \vec{\theta}^T \hat{Q} \vec{\theta} + R_n(\vec{\theta}) - \frac{1}{2}\beta_1^2\theta_1^2
					}
				 d\vec{\theta} \right| + O(n^{-3/2+4\varepsilon}) 
		 \int\limits_{\mathbb{R}^n}	
			e^{
			 - a \, \vec{\theta}^T \hat{Q} \vec{\theta}}
			 d\vec{\theta},
		\end{aligned}
	\end{equation} 
	where $\tilde{c}>0$ depends only on $a$, $b$, $d_1$, $d_2$ and $\sigma$.

Combining (\ref{x_Integral1}) and (\ref{x_Integral3}), we get that as $n\rightarrow \infty$	
\begin{equation}\label{x_Integral4}
		\begin{aligned}
			\int\limits_{U_n(n^{-1/2+\varepsilon})}	
				e^{
				i\, \vec{\beta}_n^T \vec{\theta} - a\, \vec{\theta}^T \hat{Q} \vec{\theta} + R_n(\vec{\theta})
				}&
			 d\vec{\theta} 
			 =\\
			 = 
			 \left( 1 + \frac{c^{(1)}}{n}\right)
			 \int\limits_{U_n(n^{-3/2+\varepsilon})}	
					&e^{
			 		 i (\vec{\beta}_n^T \vec{\theta}- \beta_1\theta_1)  - a\, \vec{\theta}^T \hat{Q} \vec{\theta} + R_n(\vec{\theta}) - \frac{1}{2}\beta_1^2\theta_1^2
					}
				 d\vec{\theta} + \\&+O(n^{-3/2+4\varepsilon}) 
			 \int\limits_{\mathbb{R}^n}	
			e^{
			 - a \, \vec{\theta}^T \hat{Q} \vec{\theta}}
			 d\vec{\theta}, 
		\end{aligned}	 
	\end{equation}
where $|c^{(1)}| \leq \tilde{c}\beta_1 \leq \tilde{c}b$.

We define 
\begin{equation}
	R_n^{(k)} = R_n -\frac{1}{2}\sum\limits_{j=1}\limits^{k} \beta_j^2\theta_j^2.
\end{equation}
We continue similarly to (\ref{x_Integral4})
\begin{equation}\label{x_Integral5}
		\begin{aligned}
			\int\limits_{U_n(n^{-1/2+\varepsilon})}	
				e^{
				i\left(\vec{\beta}_n^T \vec{\theta} - \sum\limits_{j=1}\limits^{k} \beta_j\theta_j\right)
				 - a\, \vec{\theta}^T \hat{Q} \vec{\theta} + R_n^{(k)}(\vec{\theta})
				}&
			 d\vec{\theta} 
			 =\\
			 = 
			 \left( 1 + \frac{c^{(k+1)}}{n}\right)
			 \int\limits_{U_n(n^{-3/2+\varepsilon})}	
					&e^{
			 		 i \left(\vec{\beta}_n^T \vec{\theta}- \sum\limits_{j=1}\limits^{k+1} \beta_j\theta_j\right)  
			 		 - a\, \vec{\theta}^T \hat{Q} \vec{\theta} + R_n^{(k+1)}(\vec{\theta}) 
					}
				 d\vec{\theta} + \\&+O(n^{-3/2+4\varepsilon}) 
			 \int\limits_{\mathbb{R}^n}	
			e^{
			 - a \, \vec{\theta}^T \hat{Q} \vec{\theta}}
			 d\vec{\theta}, 
		\end{aligned}	 
	\end{equation}
where $|c^{(k+1)}| \leq \tilde{c}\beta_k \leq \tilde{c}b$.

Combining all inequalities of (\ref{x_Integral5}) for $k = 0,1,\ldots, n-1$, we get that as $n \rightarrow \infty$
\begin{equation}\label{x_Integral6}
		\begin{aligned}
			\int\limits_{U_n(n^{-1/2+\varepsilon})}	
				e^{
				i\,\vec{\beta}_n^T \vec{\theta}
				 - a\, \vec{\theta}^T \hat{Q} \vec{\theta} + R_n(\vec{\theta})
				}
			 d\vec{\theta} 
			 &=\\
			 = 
			 \left( 1 + \frac{c^{(1)}}{n}\right)\cdots\left( 1 + \frac{c^{(n)}}{n}\right)
			 &\int\limits_{U_n(n^{-3/2+\varepsilon})}	
					e^{			 		 
			 		 - a\, \vec{\theta}^T \hat{Q} \vec{\theta} + \tilde{R}_n(\vec{\theta}) 
					}
				 d\vec{\theta} + \\&+O(n^{-1/2+4\varepsilon}) 
			 \int\limits_{\mathbb{R}^n}	
			e^{
			 - a \, \vec{\theta}^T \hat{Q} \vec{\theta}}
			 d\vec{\theta}. 
		\end{aligned}	 
	\end{equation}
Since $|c^{(k)}| \leq b\tilde{c}$ for $k = 0,1,\ldots, n-1$, using (\ref{x_Integral6}),  we obtain (\ref{eq_Lemma_x_Integral}).

\end{Proof}

\begin{Proof} {\it Lemma \ref{Lemma_Integral}.}
Note that for $\vec{\theta} \in U_n(\frac{4}{\sigma}n^{-1/2+\varepsilon})$ as $n \rightarrow \infty$
\begin{equation}
	\left\|\frac{\partial}{\partial\vec{\theta}}\sum\limits_{(v_j,v_k)\in EG} \Delta_{jk}^4\right\|_{\infty} = O(n^{-1/2+3\varepsilon}).
\end{equation}
We define $\vec{\beta}_n = Q \vec{\alpha_n}$. 
Using Lemma \ref{Lemma_x_Integral}, we find that as $n \rightarrow \infty$
\begin{equation}\label{Integral-1}
	\begin{aligned}
		\int\limits_{U_n(n^{-1/2+\varepsilon})}	
			\exp\left(
			i\,\vec{\theta}^T Q \vec{\alpha_n} - a \vec{\theta}^T \hat{Q} \vec{\theta}
			- b\sum\limits_{(v_j,v_k)\in EG} \Delta_{jk}^4 + R_n(\vec{\theta})
			\right)
		 d\vec{\theta} = \\= 
		 \Theta_{\tilde{k}_1,\tilde{k}_2} \left(
		 \int\limits_{U_n(n^{-1/2+\varepsilon})}	
			\exp\left(
			 - a \vec{\theta}^T \hat{Q} \vec{\theta}
			- b\sum\limits_{(v_j,v_k)\in EG} \Delta_{jk}^4 + \tilde{R}_n(\vec{\theta})
			\right)
		 d\vec{\theta}
		 \right)+\\
		 +O(n^{-1/2+4\varepsilon}) 
			 \int\limits_{\mathbb{R}^n}	
			e^{
			 - a \, \vec{\theta}^T \hat{Q} \vec{\theta}},
	\end{aligned}
\end{equation}
	where $\tilde{R}_n = R_n - \frac{1}{2}\sum\limits_{j=1}\limits^{n} \beta_j^2\theta_j^2$ 
	and constants $\tilde{k}_1,\tilde{k}_2$ depend only on $a$, $b$, $d_1$, $d_2$ and $\sigma$. 
	Note that for some $d_3>0$, depending only on $c$ and $\sigma$,
\begin{equation}	
	\frac{1}{2}\sum\limits_{j=1}\limits^{n}\beta_j^2\theta_j^2 \leq d_3\frac{\vec{\theta}^T \hat{Q} \vec{\theta}}{n}.
\end{equation}
	Combining (\ref{Integral_hatQ}), (\ref{Integral-1}) and Lemma \ref{Lemma_x^4_Integral}, we find that as $n\rightarrow\infty$
\begin{equation}\label{Integral0}
\begin{aligned}
	\int\limits_{U_n(n^{-1/2+\varepsilon})}	
			\exp\left(
			i\,\vec{\theta}^T Q \vec{\alpha_n} - a \vec{\theta}^T \hat{Q} \vec{\theta}
			- b\sum\limits_{(v_j,v_k)\in EG} \Delta_{jk}^4 + R_n(\vec{\theta})
			\right)
		 d\vec{\theta} = \\= 
		 \Theta_{k_1',k_2'} \left(
		  2^{\frac{n-1}{2}}\pi^{\frac{n-1}{2}} \Big/ \sqrt{\det\hat{Q}}
		 \right),
\end{aligned}
\end{equation}
where constants $k_1',k_2'$ depend only on $a$, $b$, $d_1$, $d_2$ and $\sigma$. Note that
\begin{equation}\label{Integral1}
	\begin{aligned}
		\int\limits_{-\infty}\limits^{+\infty} e^{-anx^2} dx\int\limits_{L \cap V_0}	
			\exp\left(
			i\,\vec{\theta}^T Q \vec{\alpha_n} - a\sum\limits_{(v_j,v_k)\in EG}\Delta_{jk}^2
			- b\sum\limits_{(v_j,v_k)\in EG} \Delta_{jk}^4 + R_n(\vec{\theta})
			\right)
		 dL =\\ =
		 \int\limits_{P(\theta) \in L \cap V_0}	
			\exp\left(
			i\,\vec{\theta}^T Q \vec{\alpha_n} - a \vec{\theta}^T \hat{Q} \vec{\theta}
			- b\sum\limits_{(v_j,v_k)\in EG} \Delta_{jk}^4 + R_n(\vec{\theta})
			\right)
		 d\vec{\theta}.
	\end{aligned}
	\end{equation}
	We have that 
	\begin{equation}
		\left|\exp\left(
			i\,\vec{\theta}^T Q \vec{\alpha_n} - a \vec{\theta}^T \hat{Q} \vec{\theta}
			- b\sum\limits_{(v_j,v_k)\in EG} \Delta_{jk}^4 + R_n(\vec{\theta})
			\right)\right| \leq e^{
			 - a \vec{\theta}^T \hat{Q}\vec{\theta} + \frac{d_1}{n}\vec{\theta}^T \hat{Q}\vec{\theta}
			}.	
\end{equation}
	Thus, combining (\ref{1/2U_n}), (\ref{U_n_in_P}), (\ref{Integral0}) and (\ref{Integral1}), we obtain (\ref{eq_Lemma_Integral})
\end{Proof}	
	
%%%%%%%%%%%%%%%%%%%%%%%%%%%%%%%%%%%%%%%%%%%%%%%%%%%%%%%%%%%%%%%%%%%%%%%%%%%%%%%%%%%%%%%%%%%%%%%%%%%%%%%%%%%%%%%%%%%%%%%%%
%%%%%%%%%%%%%%%%%%%%%%%%%%%%%%%%%%%%%%%%%%%%%%%%%%%%%%%%%%%%%%%%%%%%%%%%%%%%%%%%%%%%%%%%%%%%%%%%%%%%%%%%%%%%%%%%%%%%%%%%%
%%%%%%%%%%%%%%%%%%%%%%%%%%%%%%%%%%%%%%%%%%%%%%%%%%%%%%%%%%%%%%%%%%%%%%%%%%%%%%%%%%%%%%%%%%%%%%%%%%%%%%%%%%%%%%%%%%%%%%%%%
%%%%%%%%%%%%%%%%%%%%%%%%%%%%%%%%%%%%%%%%%%%%%%%%%%%%%%%%%%%%%%%%%%%%%%%%%%%%%%%%%%%%%%%%%%%%%%%%%%%%%%%%%%%%%%%%%%%%%%%%%

\section{Final remarks}

	In fact, using Lemma \ref{Lemma_using_Tutte} and Lemma \ref{Lemma_S1234}, the estimation of the number of Eulerian circuits
	 is reduced (see proof of Lemma \ref{Lemma_S0}) to estimating the integral
\begin{equation}\label{888}
		\int\limits_{V_0}	
		\exp\left(
				i\,\vec{\theta}^T Q \vec{\alpha} - \frac{1}{2}\sum\limits_{(v_j,v_k)\in EG}\Delta_{jk}^2
			- \frac{1}{12}\sum\limits_{(v_j,v_k)\in EG} \Delta_{jk}^4 + \frac{1}{2}tr(\Lambda \hat{Q}^{-1})^2
			\right) 
		 d\vec{\theta},
\end{equation}
where $\alpha$ denotes the vector composed of the diagonal elements of $\hat{Q}^{-1}$, $\Lambda$ denotes
the diagonal matrix whose diagonal elements  are equal to components of the vector $Q\vec{\theta}$. Apparently, it is possible 
to estimate integral (\ref{888})
 more accurately for particular classes of graphs and obtain  
asymptotic formulas for $Eul(G)$, similar to (\ref{Brendanf}).

Finally, we want to note that the following expression
\begin{equation}\label{my_formula}
		  2^{|EG|-\frac{n-1}{2}} \pi^{-\frac{n-1}{2}}\sqrt{t(G)} \,
		 \prod\limits_{j=1}^n \left( \frac{d_j}{2} -1 \right)!
\end{equation}
gives a surprisingly good estimate for the number of Eulerian circuits in graphs. Namely, we calculated the exact numbers 
of Eulerian circuits for small random graphs and in all cases the values given by (\ref{my_formula}) differ from the exact ones 
 within  not more than 30\% error.

\section*{Acknowledgements}
This work was carried out under the supervision of S.P. Tarasov and supported in part by RFBR grant no 11-01-00398a.
%This work was supported by the RFBR grant .

\end{document}